\documentclass{article}

\usepackage{PRIMEarxiv}
\usepackage{amsthm,amssymb}
\usepackage[utf8]{inputenc} 
\usepackage[T1]{fontenc}    
\usepackage{hyperref}       
\usepackage{url}            
\usepackage{booktabs}       
\usepackage{caption}
\usepackage{bbm}
\usepackage{comment}
\usepackage{amsfonts}       
\usepackage{nicefrac}       
\usepackage{microtype}      
\usepackage{lipsum}
\usepackage{fancyhdr}       
\usepackage{graphicx}       
\graphicspath{{media/}}     
\usepackage[fleqn]{amsmath}
\usepackage{paralist}
\usepackage[misc]{ifsym}
\usepackage{epsfig} 
\usepackage{epstopdf} 
\usepackage{xcolor}
\usepackage{tikz}
\usepackage[
backend=biber
]{biblatex}
 \usepackage{pgfplots}
     \pgfplotsset{
        table/search path={TikZData},
    }
\pgfplotsset{compat=newest}

\newtheorem{theorem}{Theorem}[section]

\newtheorem{remark}[theorem]{Remark}
\newtheorem{algorithm}[theorem]{Algorithm}

\newcommand{\R}{\mathbb{R}}
\newcommand{\N}{\mathbb{N}}
\newcommand{\de}{\mathrm{d}}

\newcommand{\e}{\mathrm{e}}
\newcommand{\A}{\boldsymbol{A}}

\newcommand{\Height}{6.5 cm}
\newcommand{\Width}{6.5 cm}
\newcommand{\Heightgs}{3.8 cm}
\newcommand{\Widthgs}{3.8 cm}
\newcommand{\Heightt}{5.8 cm}
\newcommand{\Widtht}{5.8 cm}

\definecolor{plotcolor1}{RGB}{228,26,28}    
\definecolor{plotcolor2}{RGB}{55,126,184}   
\definecolor{plotcolor3}{RGB}{77,175,74}    
\definecolor{plotcolor4}{RGB}{152,78,163}   
\definecolor{plotcolor5}{RGB}{255,127,0}    
\definecolor{plotcolor6}{RGB}{255,255,51}   

\usepackage{cases}   
\title{Comparing intrusive and non-intrusive polynomial chaos for a class of exponential time differencing schemes
}

\author{
  Julian Clausnitzer \\
  J\"{u}lich Supercomputing Centre, Forschungszentrum J\"{u}lich GmbH, 52425 J\"{u}lich \\
  {j.clausnitzer@fz-juelich.de} \\
   \AND
  Andreas Kleefeld \\
    J\"{u}lich Supercomputing Centre, Forschungszentrum J\"{u}lich GmbH, 52425 J\"{u}lich\\
  Faculty of Medical Engineering and Technomathematics, University of Applied Sciences Aachen, 52428 J\"{u}lich\\
  {a.kleefeld@fz-juelich.de} \\
}

\begin{document}

\definecolor{plotcolor1}{RGB}{120,94,240}
\definecolor{plotcolor2}{RGB}{100,143,255}
\definecolor{plotcolor3}{RGB}{206,142,0}
\definecolor{plotcolor4}{RGB}{254,97,0}
\definecolor{plotcolor5}{RGB}{220,38,127}
\definecolor{plotcolor6}{RGB}{0,38,127}

\pgfplotscreateplotcyclelist{goodcolors}{
{color=plotcolor1, style=solid, line width=2pt},
{color=plotcolor2, style=dashed, line width=2pt},
{color=plotcolor3, style=dotted, line width=2pt},
{color=plotcolor4, style=densely dashed, dash pattern=on 4pt off 2pt, line width=2pt},
{color=plotcolor5, style=dashed, dash pattern=on 6pt off 2pt, line width=2pt},
{color=plotcolor6, style=solid, line width=1pt},
}

\pgfplotscreateplotcyclelist{goodcolorsalt}{
{color=plotcolor2, style=dashed, line width=2pt},
{color=plotcolor3, style=dotted, line width=2pt},
{color=plotcolor5, style=dashed, dash pattern=on 6pt off 2pt, line width=2pt},
{color=plotcolor6, style=solid, line width=1pt},
}

\maketitle



\begin{abstract}
We consider the numerical approximation of different ordinary differential equations (ODEs) and partial differential equations (PDEs) with periodic boundary conditions involving a one-dimensional random parameter, comparing the intrusive and non-intrusive polynomial chaos expansion (PCE) method. We demonstrate how to modify two schemes for intrusive PCE (iPCE) which are highly efficient in solving nonlinear reaction-diffusion equations: A second-order exponential time differencing scheme (ETD-RDP-IF) as well as a spectral exponential time differencing fourth-order Runge-Kutta scheme (ETDRK4). In numerical experiments, we show that these schemes show superior accuracy to simpler schemes such as the EE scheme for a range of model equations and we investigate whether they are competitive with non-intrusive PCE (niPCE) methods. We observe that the iPCE schemes are competitive with niPCE for some model equations, but that iPCE breaks down for complex pattern formation models such as the Gray-Scott system.
\end{abstract}
 \section{Introduction.}\label{section:introduction}
In this work, we aim to solve time-dependent random differential equations of the form
\begin{align}\label{differentialequation}
\begin{cases}
    \frac{\partial u(x,t,\omega)}{\partial t} &= D\Delta u(x,t,\omega) + F(\omega,u(x,t,\omega)),\\
    u(x,0,\omega) &= u_\mathrm{init}(x),
\end{cases}
\end{align}
with periodic boundary conditions, where $t\in (0,T]$ for some $T>0$, $x\in (-1,1)^d$ for $d\in \{1,2\}$, and $\omega\in\Omega$ for a probability space $(\Omega,\mathcal{F},\mathbb{P})$. The constant $D$ may be positive, so that we deal with a PDE, or $D=0$, in which case \eqref{differentialequation} is an ODE. The potentially nonlinear function $F$ depends on $\omega$, whereas the initial condition $u_\mathrm{init}$ does not, i.e. we will be working with a deterministic initial condition and a random nonlinearity.\\
 Using PCE representations to approximate random processes has been investigated in stochastic mechanics \cite{ghanemspanos1} and computational fluid dynamics (CFD) \cite{debusschere2003,maitremultigrid,lucor2004} for problems involving random parameters. With the usage of Wiener-Askey polynomial chaos \cite{xiukarniadakis}, it became possible to not only efficiently treat Gaussian inputs, but other input distributions as well. However, especially intrusive PCE (iPCE), which relies on inserting the PCE into the equation such as \eqref{differentialequation} and subsequently discretizing the resulting system, has been facing issues with long-term integration as well as sharp dependencies in the stochastic parameter space (see \cite[p. 1]{eckert2020polynomial} and \cite[pp. 1]{kaur2022adaptive}). This has also been noted for unsteady flow problems in CFD \cite{wan2006long,gerritsma2010time,le2010asynchronous}. Since numerical solvers for deterministic problems have to be modified considerably in order to be used for iPCE, a relevant question is to what extent the benefits of these solvers carry over to the iPCE case. We consider a few model equations with a diffusion term and a nonlinear term, and we also look at a reaction-diffusion equation, the Gray-Scott model. While there are many numerical methods to solve reaction-diffusion equations (see \cite{chiu1997adi,shakeri2011finite,song2021convergence,zhang2008second}), in this paper we consider the following two numerical schemes: Firstly, a second-order finite difference exponential time differencing scheme which approximates matrix exponentials using rational functions having real distinct poles combined with an integrating factor approach (ETD-RDP-IF) \cite{Kleefeld2020}. Secondly, we look at a spectral exponential time differencing Runge-Kutta fourth order scheme (ETDRK4) \cite{kassamsolving}. The performances of these schemes are compared along with a simple explicit Euler (EE) scheme. We show that the iPCE schemes based on ETD-RDP-IF and ETDRK4 are competitive with their niPCE counterparts, but that more complex models such as the Gray-Scott model cannot be treated with the iPCE implementation presented here.\\
 The contribution of this paper is a development and implementation of the iPCE method using the ETD-RDP-IF and ETDRK4 algorithms and a thorough quantitative investigation of the errors. The iPCE errors are compared to errors produced by Monte Carlo simulations and Gaussian Quadrature which represent niPCE simulations.
 The MATLAB program code to generate all images and error plots in this paper can be found at \url{https://github.com/JulianSPDE/PCE}. The error plots can be reproduced by executing \texttt{runall.m}; more detailed instructions can be found in a readme file.\\
The paper is organized as follows. In Chapter \ref{section:pce}, we describe the mathematical setup and give a short introduction into non-intrusive and intrusive polynomial chaos expansion. The quantities of interest are then the expectation and the variance of the time-dependent stochastic partial/ordinary differential equation. Chapter \ref{section:numerical_schemes} reviews different numerical schemes to solve these deterministic equations such as the explicit Euler, the second-order exponential time differencing method with distinct real poles, and the spectral exponential time differencing fourth-order Runge-Kutta method. We then develop in detail the corresponding methods for the non-deterministic equations. Extensive numerical results are given in Chapter \ref{section:numericalexperiments}. First, results are presented for a one-dimensional non-deterministic ordinary and partial differential equation involving a linear, a quadratic, and a cubic term. Second, results are given for a two-dimensional equation as well as for a system of two-dimensional equations --- the Gray-Scott model. In Chapter \ref{section:Outlook}, a conclusion and an outlook is given.
\section{Polynomial chaos expansion.}\label{section:pce}
Let $(\Omega,\mathcal{F},\mathbb{P})$ be a probability space with state space $\Omega$, a $\sigma$-algebra $\mathcal{F}$ and a probability measure $\mathbb{P}$, and let $(H,\mathcal{H})$ be a measurable space where $H$ is a Hilbert space and $\mathcal{H} = \mathcal{B}(H)$ is the Borel $\sigma$-algebra on $H$. Suppose that $X$ is in the space $L^2(\Omega,\mathcal{F},\mathbb{P};H)$, the space of square-integrable $H$-valued random variables $X:\Omega\rightarrow H$, i.e. $\mathbb{E}[\|X\|_H^2]<\infty.$ Suppose also that $\xi\in L^2(\Omega,\mathcal{F},\mathbb{P};\R)$ is a real-valued random variable and $X$ depends on $\xi$ by some function $g$, i.e. $X(\omega) = g(\xi(\omega)).$ Using the notation from \eqref{differentialequation}, for a given fixed time $t>0$ we understand $X(t) = u(x,t,\omega)$ as a time-dependent element of $L^2(\Omega,\mathcal{F},\mathbb{P};H)$. Then, the generalized polynomial chaos expansion (gPCE) for $u(x,t,\omega)$ is given by \cite[Eq. (3.1)]{xiukarniadakis}
\begin{align}\label{pce}
    u(x,t,\omega) = \sum_{i=0}^\infty  u_i(x,t) P_i(\xi(\omega)),
\end{align}
where $u_i\in H$ are time-dependent, deterministic coefficient functions which have to be determined, and $\{P_i\}_{i\in\mathbb{N}_0}$ is a family of orthogonal polynomials. From now, we will usually drop the dependency on $\omega$ and instead just write $\xi$ or leave $\omega$ and $\xi$ entirely when it is clear from the context that we deal with a random quantity. According to Wiener-Askey polynomial chaos, different distributions for $\xi$ lead to different families $\{P_i\}_{i\in\mathbb{N}_0}$, of which we recall two in Table \ref{orthopoly}. It is known that \eqref{pce} converges in the $L^2$ sense \cite[pp. 101--119]{ghanemspanos1}, and since the polynomials $\{H_i\}_{i\in\N_0}$ are orthogonal, the coefficient functions $u_i$ can be determined by the orthogonal $L^2$ projection
\begin{align}\label{projection}
    u_i(x,t) = \frac{\mathbb{E}[u(x,t,\cdot)P_i(\xi)]}{\mathbb{E}[P_i(\xi)^2]}.
\end{align}
We will assume that the orthogonal polynomials are normed, so that $\langle P_i^2\rangle = \mathbb{E}[P_i(\xi)^2] = 1$, $i=0,\dots, N$. We note that in this paper, subscripts will usually refer to the index of the PCE coefficients, while superscripts usually refer to the time step, so for example $u_i^n$ refers to the $i$-th PCE coefficient at time $t_n$.\\
In order to solve equation \eqref{differentialequation}, the gPCE \eqref{pce} can be utilized in two different ways:
\begin{enumerate}
    \item Non-intrusive PCE is recalled in Section \ref{section:nonintrusivepce}. It relies on sampling the input random variable $\xi$, yielding samples $\xi_1,\dots,\xi_M$. The deterministic equation
\begin{align*}
\begin{cases}
    \frac{\partial u(x,t,\xi_j)}{\partial t} &= D\Delta u(x,t,\xi_j) + F(\xi_j,u(x,t,\xi_j)),\\
    u(x,0,\xi_j) &= u_\mathrm{init}(x),
\end{cases}
\end{align*}
is solved for $j=1,\dots,M$ using deterministic solvers which are treated as a `black box'. The samples $u(x,t,\xi_j)$ are then used to numerically compute the expectations in \eqref{projection}.
\item Intrusive PCE is treated in Section \ref{section:intrusivepce}. It plugs the gPCE \eqref{pce} into \eqref{differentialequation} in the beginning, yielding
\begin{align*}
    \begin{cases}
        \frac{\partial}{\partial t} \sum_{i=0}^\infty u_i(x,t)P_i(\xi) = D\Delta \sum_{i=0}^\infty u_i(x,t)P_i(\xi) + F\left(\omega,\sum_{i=0}^\infty u_i(x,t)P_i(\xi)\right) \\
        u_0(x,0) = u_\mathrm{init}(x),\quad u_\eta(x,0)=0,\ \eta >0
    \end{cases}\\
    \Leftrightarrow \begin{cases}
         \sum_{i=0}^\infty \frac{\partial}{\partial t}u_i(x,t)P_i(\xi) = D \sum_{i=0}^\infty \Delta u_i(x,t)P_i(\xi) + F\left(\omega,\sum_{i=0}^\infty u_i(x,t)P_i(\xi)\right) \\
        u_0(x,0) = u_\mathrm{init}(x),\quad u_\eta(x,0)=0,\ \eta>0.
    \end{cases} 
\end{align*}
After truncating the series to $N+1$ functions (index $0$ to $N$) and performing a Galerkin projection, i.e. for a fixed index $\eta\in\N_0$ multiplying with $P_\eta(\xi)$ and taking $\mathbb{E[\cdot]}$, using the orthogonality property $\mathbb{E}[P_i(\xi)P_j(\xi)] = \delta_{ij}$, we obtain
\begin{align}\label{aftergalerkin}
    \begin{cases}
         \frac{\partial}{\partial t}u_\eta(x,t)= D \Delta u_\eta(x,t) + \mathbb{E}\left[F\left(\omega,\sum_{i=0}^N u_i(x,t)P_i(\xi)\right)P_\eta(\xi)\right],\quad \eta=0,\dots, N \\
        u_0(x,0) = u_\mathrm{init}(x),\quad u_i(x,0)=0,\ i>0.
    \end{cases}
\end{align}
The system \eqref{aftergalerkin} is then discretized and is solved using a solver specifically implemented for the iPCE system. Its implementation depends heavily on the chosen discretization in space and time, and on the function $F$. We will treat a finite difference discretization of second order in space with an explicit Euler scheme which is first-order in time in Section \ref{section:expliciteuler} and with a second-order exponential time differencing scheme in Section \ref{section:etdrdpif}, as well as a spectral discretization in space with a fourth-order Runge-Kutta time stepping scheme in Section \ref{section:etdrk4}.
\end{enumerate}
We will mainly be interested in computing the expected value $\mathbb{E}[u(x,t,\cdot)]$ and the variance $\mathrm{Var}[u(x,t,\cdot)]$ of the solution. We note that due to the properties of the orthogonal polynomials $\{P_i\}_{i\in\mathbb{N}_0}$, it is
\begin{align}\label{meanandvariance}
    \mathbb{E}[u(x,t,\cdot)] = u_0(x,t),\qquad \mathrm{Var}[u(x,t,\cdot)] = \sum_{i=1}^\infty |u_i(x,t)|^2.
\end{align}
We also note that the periodic boundary conditions for $u$ are transferred directly to the PCE base functions: Inserting the periodic boundary condition $u(-1,t)=u(1,t)$, $t\in [0,T]$, into \eqref{pce} and then multiplying with $P_j(\xi)$ and taking $\mathbb{E}[\cdot]$ yields
\begin{align}\label{boundarycondition_pce}
    \sum_{i=0}^\infty u_i(-1,t)P_i(\xi) = \sum_{i=0}^\infty u_i(1,t)P_i(\xi)\ \Rightarrow\ u_j(-1,t)=u_j(1,t),\quad j\in\N_0,\ t\in [0,T]
\end{align}
due to the orthonormality of $\{P_i\}_{i\in \N_0}$. In higher spatial dimensions, the above argument can easily be replicated by applying \eqref{boundarycondition_pce} to any pair of points $x_1$, $x_2$ on the boundary for which $u(x_1,t) = u(x_2,t)$ holds.
\begin{table}
    \centering
    \begin{tabular}{|c|c|c|c|}
    \toprule
         Distribution of $\xi$ & Polynomial $P_n$ & Weight function $\rho(x)$ & Support $S$ \\
         \midrule
         Gaussian & Hermite $H_n$ & $\e^{-x^2}$ & $(-\infty,\infty)$ \\
         Uniform & Legendre $L_n$ & $\frac{1}{2}$ & $[-1,1]$ \\
         \bottomrule
    \end{tabular}
    \caption{Continuous probability distributions associated with families of orthogonal polynomials, table from \cite[Table 4.1]{xiukarniadakis}. Note that Legendre polynomials can be linearly rescaled so that they are supported on an arbitrary closed interval $[a,b]\subset \R$. We will later exclusively use Legendre polynomials supported on $[a,b]\subset \R$.}
    \label{orthopoly}
\end{table}
\section{Numerical schemes.}\label{section:numerical_schemes}
 Throughout this paper, we will be using $p$ points in each spatial dimension, $M$ time steps, a time step size $k=T/M$, a spatial step size $h=p^{-1}$ and $N+1$ PCE coefficients (with the index ranging from 0 to $N$). For the random equations \eqref{differentialequation}, we will assume that the random input is uniformly distributed, so we will always use Legendre polynomials in the gPCEs. The cases for other random distributions and the corresponding family of orthogonal polynomials according to Table \ref{orthopoly} work in an analogous fashion.
 \subsection{The schemes for deterministic PDEs}
 We will now recall the deterministic numerical schemes for solving the deterministic initial value problem
 \begin{align*}
     \begin{cases}
         \frac{\partial u(x,t)}{\partial t} = D\Delta u(x,t) + F(u(x,t)),\\
         u(x,0) = u_{\mathrm{init}}(x).
     \end{cases}
 \end{align*}
 These schemes can be used directly for non-intrusive PCE (see Section \ref{section:nonintrusivepce}) or modified to be used in intrusive PCE (see Section \ref{section:intrusivepce}). We start with the well-known explicit Euler scheme in Section \ref{section:expliciteuler}, for example as introduced in \cite[pp. 240]{hoffman2018numerical}, secondly recall the ETD-RDP and ETD-RDP-IF schemes in Section \ref{section:etdrdpif} presented in \cite{Kleefeld2020}, and finally the ETDRK4 scheme presented in \cite{kassamsolving}.
 \subsubsection{Explicit Euler scheme}\label{section:expliciteuler}
 For the explicit Euler scheme and the ETD-RDP-IF scheme, we will use a second-order central finite difference discretization in space. We denote by $p$ the spatial resolution in each dimension, so that in one dimension we have the discretized function $\boldsymbol{u}(t) = (u(x_0,t),u(x_1,t),\dots,u(x_p,t))^\top$ where $x_i = -1+2i/p$, $i=0,\dots, p$, $t\in[0,T]$. The second-order finite difference approximation of the Laplacian with periodic boundary conditions is given by the matrix
 \begin{align*}
     A_p = -\frac{1}{h^2}\begin{pmatrix}
         -2 & 1 & & & 1 \\
         1 & -2 & 1 & & \\
         & \ddots & \ddots & \ddots & \\
         & & 1 & -2 & 1 \\
         1 & & & 1 & -2
     \end{pmatrix}\in\R^{p\times p}.
 \end{align*}
Also, the discretization of the two-dimensional Laplacian is given by $A := I_{N+1}\otimes A_p + A_p \otimes I_{N+1}$. Denoting $\boldsymbol{u}^n = (u(x_0,t_n),u(x_1,t_n),\dots,u(x_p,t_n))^\top$ for $n=0,\dots,M$, $t_n=n\cdot k$, the explicit Euler scheme for solving \eqref{differentialequation} is then given by
 \begin{algorithm}[Deterministic explicit Euler scheme]\label{algorithm:EE_deterministic}
 \begin{align*}
     \boldsymbol{u}^{n+1} = \boldsymbol{u}^n + k(DA_p\boldsymbol{u}^n + F(\boldsymbol{u}^n)),\quad n=0,\dots, M-1.
 \end{align*}    
 \end{algorithm}
 \subsubsection{ETD-RDP-IF}\label{section:etdrdpif}
 The ETD-RDP-IF scheme presented in \cite{Kleefeld2020} is a second-order exponential time-differencing (ETD) scheme which makes use of the approximation of the matrix exponentials by rational functions having real distinct poles (RDP). Additionally, if the spatial dimension is greater than one, the scheme uses dimensional splitting with an integrating factor (IF) approach, making it the ETD-RDP-IF scheme. For one spatial dimension, the ETD-RDP scheme is given by \cite[Eq. (17)]{Kleefeld2020}
 \begin{algorithm}[Deterministic ETD-RDP scheme]\label{algorithm:ETDRDP_deterministic}
 \begin{align}\label{etdrdpscheme}
     \boldsymbol{u}^{n+1} &= \left( I_p + \frac{k}{3} DA_p \right)^{-1}[9\boldsymbol{u}^n+2kF(\boldsymbol{u}^n)+kF(\boldsymbol{u}^{n+1}_*)] \\
     &- \left( I_p + \frac{k}{4}DA_p \right)^{-1}\left[ 8\boldsymbol{u}^n + \frac{3k}{2}F(\boldsymbol{u}^n) + \frac{k}{2}F(\boldsymbol{u}^{n+1}_*) \right],\\
     \boldsymbol{u}^{n+1}_* &= (I_p+kDA_p)^{-1} (\boldsymbol{u}^n + kF(\boldsymbol{u}^n)),\quad n=0,\dots,M-1,
 \end{align}    
 \end{algorithm} 
 where $I_p$ denotes the $p\times p$ identity matrix. \\
 For two spatial dimensions, the ETD-RDP-IF scheme is given by (for the full derivation of the scheme, we refer to \cite{Kleefeld2020})
 \begin{algorithm}[Deterministic ETD-RDP-IF scheme]\label{algorithm:ETDRDPIF_deterministic}
 \begin{align*}
     \boldsymbol{u}^{n+1} &= \left( I_{p^2} + \frac{k}{3} D A_2 \right)^{-1}\left[  \left\{9\left( I_{p^2} + \frac{k}{3} DA_1\right)^{-1}   - 8\left( I_{p^2} + \frac{k}{4}DA_1\right)^{-1} \right\}\right. \\ &
     \cdot \left. \left\{ 9\boldsymbol{u}^n+2kF(\boldsymbol{u}^n)\right\}+kF(\boldsymbol{u}^{n+1}_*) \right] \\
     &- \left( I_{p^2} + \frac{k}{4} D A_2 \right)^{-1}\left[  \left\{9\left( I_{p^2} + \frac{k}{3} DA_1\right)^{-1}   - 8\left( I_{p^2} + \frac{k}{4}DA_1\right)^{-1}\right\} \right.\\
     &\cdot \left. \left\{8\boldsymbol{u}^n+\frac{3k}{2}F(\boldsymbol{u}^n)\right\}+\frac{k}{2}F(\boldsymbol{u}^{n+1}_*) \right], \\
     \boldsymbol{u}^{n+1}_* &= (I_{p^2}+kDA_{p^2})^{-1}(I_{p^2}+kDA_1)^{-1} (\boldsymbol{u}^n + kF(\boldsymbol{u}^n)),\quad n=0,\dots,M-1,
 \end{align*}
 \end{algorithm}
 where $A_1 := I_p\otimes A_p$, $A_2 := A_p\otimes I_p$.
 \subsubsection{ETDRK4}\label{section:etdrk4}
 The ETDRK4 scheme has been shown in \cite{kassamsolving} to be a powerful method for solving reaction-diffusion equations which is superior to standard finite difference schemes. For its spatial discretization, it uses a spectral approach, making use of the discrete Fourier transform (DFT) which can be implemented using the well-known fast Fourier transform (FFT). Using the shorthand $u(x_n)=:\boldsymbol{u}_n$, $n=0,\dots, p$, the DFT and inverse DFT formulae are given by
 \begin{align*}
     \widehat{ \boldsymbol{u}}_j = \sum_{n=0}^{p-1} \e^{-2\pi\mathrm{i} j x_n} \boldsymbol{u}_n,\quad j=-\frac{p}{2}+1,\dots, \frac{p}{2},\quad \boldsymbol{u}_n = \frac{1}{p}\sum_{j=-p/2+1}^{p/2}\e^{2\pi\mathrm{i}jx_n}\widehat{\boldsymbol{u}}_j,\quad n=0,\dots, p-1.
 \end{align*}
We will later also denote $\hat{\boldsymbol{u}} =: \mathcal{F}(\boldsymbol{u})$, especially when emphasizing that a transformation from physical space to Fourier space is performed, rather than a mere manipulation of terms in Fourier space. The derivatives on the spatial grid can now be approximated conveniently by noting that differentiation $\boldsymbol{u}' = (u'(x_n))_{n=0,\dots,p}$ in the physical space corresponds to a simple multiplication in Fourier space: $\widehat{\boldsymbol{u}'} = (\mathrm{i}j\hat{u}_j)_{j=-\frac{p}{2}+1,\dots,\frac{p}{2}}$. In particular, the Laplace operator in Fourier space ends up being diagonal: $\widehat{\Delta \boldsymbol{u}} = L\cdot \hat{\boldsymbol{u}}$ for a diagonal matrix $L$ with $L_{jj} = -j^2$, $j=-\frac{p}{2}+1,\dots,\frac{p}{2}$. For the ETDRK4 scheme, the DFT is applied to the discretized function $\boldsymbol{u}$, the time stepping using the fourth-order Runge-Kutta method is carried out in Fourier space and after the last time step the solution is transformed back to physical space. The formulae for the Runge-Kutta fourth order scheme can be found in \cite[p. 3]{kassamsolving}. One subtlety to note is that in each time step, for each evaluation using the nonlinear function $F$, the argument $\widehat{\boldsymbol{u}^n}$ needs to be transformed to physical space, then $F$ is evaluated, and $F(\mathcal{F}^{-1}(\widehat{\boldsymbol{u}^n}))$ is transformed back into Fourier space (see the scheme below). Taking this into account, we introduce the shorthand $\hat F(\widehat{\boldsymbol{u}^n}) := \mathcal{F}(F(\mathcal{F}^{-1}(\widehat{\boldsymbol{u}^n})))$ and the full ETDRK4 scheme is given by
\begin{algorithm}[Deterministic ETDRK4 scheme]\label{algorithm:ETDRK4_deterministic}
\begin{align}\label{etdrk4_deterministic}
    \widehat{\boldsymbol{u}^0} &= \mathcal{F}(\boldsymbol{u}^0),\\ \nonumber
    \widehat{\boldsymbol{u}^{n+1}} &= \e^{Lh}\widehat{\boldsymbol{u}^n} + h^{-2}L^{-3}\{[-4\cdot I_p-Lh+\e^{Lh}(4\cdot I_p-3Lh+(Lh)^2)]\hat F(\widehat{\boldsymbol{u}^n}) \\ \nonumber
    &+ 2[2\cdot I_p+Lh+\e^{Lh}(-2\cdot I_p+Lh)](\hat F(a_n)+\hat F(b_n) \\ \nonumber
    &+ [-4-3Lh-(Lh)^2+\e^{Lh}(4\cdot I_p-Lh)]\hat F(c_n)\},\\ \label{deterministic_an}
    a_n &= \e^{Lh/2}\widehat{\boldsymbol{u}^n} + L^{-1}(\e^{Lh/2}-I_p)\hat F(\widehat{\boldsymbol{u}^n})), \\ \label{deterministic_bn}
    b_n &= \e^{Lh/2}\widehat{\boldsymbol{u}^n} + L^{-1}(\e^{Lh/2}-I_p)\hat F(a_n),\\\label{deterministic_cn}
    c_n &= \e^{Lh/2}a_n + L^{-1}(\e^{Lh/2}-I_p)(2\hat F(b_n)-\hat F(\widehat{\boldsymbol{u}^n})),\quad n=0,\dots,M-1,\\ \nonumber
    \boldsymbol{u}^M &= \mathcal{F}^{-1}(\widehat{\boldsymbol{u}^M}).
\end{align}
\end{algorithm}
We note that the scheme presented in this form was introduced by Cox and Matthews \cite{cox2002exponential}. However, it can lead to numerical instabilities resulting from cancellation errors in the expressions of $a_n$, $b_n$ and $c_n$ \cite[p. 6]{kassam2005fourth}. Instead of evaluating the critical expressions directly, Cauchy's integral formula
\begin{align}\label{integralformula}
    f(L) = \frac{1}{2\pi \mathrm{i}} \int_\Gamma f(t) (tI_p-L)^{-1}\; \de t
\end{align}
is used, where $\Gamma$ is a contour which encloses the eigenvalues of $L$. The trapezoidal rule is suitable to evaluate this integral, since it converges exponentially for complex contour integrals \cite{davis1959numerical}. In our case, $L$ is diagonal and the contours for \eqref{integralformula} to evaluate $a_n$, $b_n$ and $c_n$ may simply be chosen elementwise, so that for each diagonal element $L_{ii}$ of $L$, $i=1,\dots, p$, we pick one circle around $L_{ii}$. The approximated integral \eqref{integralformula} then simplifies to
\begin{align}\label{complexrootsmean}
    \frac{1}{R}\sum_{i=1}^R f(L_{ii}+r_i),
\end{align}
where $r_1,\dots, r_R$ are the complex roots of unity of order $R$ sitting on the unit circle shifted by $L_{ii}$.\\
Furthermore, in order to avoid errors caused by aliasing in Fourier space, anti-aliasing is needed in the program code. We refer to \cite[p. 11]{kassamsolving} on how to do this appropriately.
 \subsection{Non-intrusive PCE}\label{section:nonintrusivepce}
 The idea of non-intrusive PCE is to compute the coefficient functions $u_i$ in \eqref{pce} as numerical integrals
 \begin{align}\label{pcenumint}
     u_i(x,t) = \frac{1}{\langle P_i^2 \rangle} \int_a^b u(x,t,\xi)P_i(\xi)\rho(\xi)\; \de \xi \approx \sum_{j=1}^q w_j u(x,t,\xi_j) P_i(\xi_j) \rho(\xi_j), 
     \end{align}
     (with $\langle P_i^2 \rangle = 1$, see p. 2) where the integral in \eqref{pcenumint} is a Banach space-valued integral (for example a Bochner integral, see \cite[Appendix C]{engelnagel}) and the points $\xi_j$ and weights $w_j$, $j=1,\dots, q$ are chosen according to the used method. We will compare three different methods:
     \begin{itemize}
         \item For the classical Monte Carlo (MC) method, $w_j=q^{-1}$ and $\xi_j$ are sampled randomly from the distribution of the random variable $\xi$. 
         \item For a quasi-Monte Carlo (QMC) method, $w_j=q^{-1}$ and $\xi_j$ are drawn from low-discrepancy sequences such as a Sobol sequence (see \cite{sobol}) or Halton sequences (see \cite{halton}). These sequences cover the parameter space $[a,b]$ more evenly and generally achieve a convergence order of almost one instead of $1/2$ for classical MC (see \cite{halton1960efficiency} and \cite[p. 980]{niederreiter}).
         \item For Gaussian quadrature (GQ), the $\xi_j$ are quadrature nodes given by the roots of $P_q$ and the quadrature weights are given by
         \begin{align*}
             w_j = \frac{a_q}{a_{q-1}} \frac{\langle P_{q-1},P_{q-1} \rangle}{P_q'(x_j)P_{q-1}(x_j)},\quad j=1,\dots, q,
         \end{align*}
         where $a_k$ is the coefficient of $x^k$ in $P_q$. An efficient and numerically stable way to compute the quadrature nodes and weights is given by the Golub-Welsch algorithm (see \cite{golub} and \cite[p. 188]{numericalrecipes}). 
     \end{itemize}
\begin{remark}
  For a multi-dimensional random input $\boldsymbol{\xi} = (\xi_1,\dots, \xi_P)$, the samples $\boldsymbol{\xi}_j = (\xi_{1,j},\dots,\xi_{P,j})$, $j=1,\dots,q$ are created independently in each component for MC, as a $P$-dimensional Sobol sequence for QMC, or as $P$-dimensional Gaussian quadrature. We discuss multidimensional inputs a bit more in Remark \ref{multidimensionalremark}. 
\end{remark}
   If we want to compute the mean $\mathbb{E}[u(x,t,\cdot)]$, we only need the first coefficient function $u_0(x,t)$, if we compute the variance we need to compute $u_1,\dots, u_N$ (see \eqref{meanandvariance}) and introduce a truncation error. In order to compute the variance, it is best to compute the samples $u(x,t,\xi_j)$, $j=1,\dots,q$, and reuse them to compute all the integrals in \eqref{pcenumint} for $i=1,\dots, q$.
 \subsection{Intrusive PCE}\label{section:intrusivepce}
 We have already given an introduction to intrusive PCE in Section \ref{section:pce}. The aim is now to develop the intrusive PCE numerical schemes needed to solve \eqref{aftergalerkin}. We will consider in this section the nonlinear function $F(\omega,u(x,t,\omega))=K(\omega)u(x,t,\omega)^3$ with a uniformly distributed constant $K(\omega) = \xi(\omega)\sim \mathcal{U}[a,b]$. We stress here that we will henceforth use $K$ and $\xi$ interchangably, and the orthogonal polynomials will be the Legendre polynomials supported on $[a,b]\subset\R$. Among others, this nonlinear function $F$ will also be treated in Section \ref{section:numericalexperiments}. We describe the explicit Euler scheme in \ref{section:expliciteuler}, the ETD-RDP scheme in Section \ref{section:intrusive_ETDRDPIF} and the ETDRK4 scheme in Section \ref{section:intrusive_etdrk4}.\\
 For both schemes, it is a necessary prerequisite to discuss the gPCE for the term $u^3$ given the PCE of $u$. Given the gPCE \eqref{pce}, we have
 \begin{align*}
     u^3(x,t,\omega) &= \left( \sum_{i=0}^\infty u_i(x,t)P_i(\xi)\right)\cdot \left( \sum_{j=0}^\infty u_j(x,t)P_j(\xi)\right)\cdot \left( \sum_{k=0}^\infty u_k(x,t)P_k(\xi)\right) \\
     &= \sum_{i=0}^\infty \sum_{j=0}^\infty \sum_{k=0}^\infty u_i(x,t)u_j(x,t)u_k(x,t)P_i(\xi)P_j(\xi)P_k(\xi).
 \end{align*}
 After a Galerkin projection (multiplying with $P_\eta$ for $\eta\in\N_0$ and taking $\mathbb{E}[\cdot]$) and making an approximation by truncating the series to $N+1$ terms each, plugging the gPCE into \eqref{aftergalerkin} we obtain for $\eta=0,\dots, N$ (see also \cite[pp. 4]{debusschere2004})
 \begin{align*}
     \begin{cases}
         \frac{\partial}{\partial t} u_\eta(x,t) = D\Delta u_\eta(x,t) + \mathbb{E}\left[K \sum_{i,j,k=0}^N u_i(x,t)u_j(x,t)u_k(x,t)P_i(\xi) P_j(\xi)P_k(\xi) P_\eta(\xi) \right]\\
         u_0(x,0) = u_\mathrm{init}(x),\quad u_\eta(x,0)=0,\ \eta>0.
     \end{cases}
 \end{align*}
 Extracting the deterministic coefficient functions, we have
 \begin{align*}
     \begin{cases}
         \frac{\partial}{\partial t} u_\eta(x,t) = D\Delta u_\eta(x,t) +  \sum_{i,j,k=0}^N u_i(x,t)u_j(x,t)u_k(x,t)\mathbb{E}\left[K P_i(\xi) P_j(\xi)P_k(\xi) P_\eta(\xi) \right]\\
         u_0(x,0) = u_\mathrm{init}(x),\quad u_\eta(x,0)=0,\ \eta>0.
     \end{cases}
 \end{align*}
 and, using from now on the notation $\boldsymbol{K}_{ijk\eta}:=\mathbb{E}[KP_i(\xi)P_j(\xi)P_k(\xi)P_\eta(\xi)]$, it is
 \begin{align}\label{intrusivepcesystem}
     \begin{cases}
         \frac{\partial}{\partial t} u_\eta(x,t) = D\Delta u_\eta(x,t) +  \sum_{i,j,k=0}^N \boldsymbol{K}_{ijk\eta} u_i(x,t)u_j(x,t)u_k(x,t)\\
         u_0(x,0) = u_\mathrm{init}(x),\quad u_\eta(x,0)=0,\ \eta>0.
     \end{cases}
 \end{align}
 The tensor $\boldsymbol{K}$ will be computed as a pre-processing step and used throughout a simulation. 
\begin{remark}\label{implementationremark}
    In our implementations, we use an equivalent way of computing the gPCE of $u^3$: We make use of an addition theorem for Legendre polynomials \cite{legendreproduct} stating that for $\alpha,\beta\in\mathbb{N}$ and $\alpha\wedge \beta = \min(\alpha,\beta)$,
    \begin{align}\label{legendreadditiontheorem}
        P_\alpha(x)P_\beta(x) &= \sum_{i\leq \alpha\wedge \beta} C(\alpha,\beta,p)P_{\alpha+\beta-2p}(x),\quad \text{where}\\
        C(\alpha,\beta,p) &= \frac{A_p A_{\alpha-p}A_{\beta-p}}{A_{\alpha+\beta-p}}\cdot \frac{2\alpha + 2\beta -4p + 1}{2\alpha + 2\beta -2p + 1}\cdot \sqrt{\frac{(2\alpha + 1)(2\beta + 1)}{2(\alpha+\beta-2p)+1}}, \\ A_r &:= \frac{(1/2)_r}{r!},\quad (a)_r := a\cdot (a+1)\cdot \dots \cdot (a+r-1),\quad (a)_0 = 1.
    \end{align}
    Using \eqref{legendreadditiontheorem} to expand the gPCE of $u^3$ yields, after applying a Cauchy product twice (and dropping $x$, $t$ and $\omega$), 
    \begin{align*}
        u^3 = \sum_{\ell=0}^\infty \sum_{m=0}^\ell \sum_{j=0}^m \sum_{p=0}^{j\wedge (m-j)} \sum_{n=0}^{(\ell-m)\wedge (m-2p)} u_{\ell-m} u_j u_{m-j} C(j,m-j,p)C(\ell-m,m-2p,n)  P_{\ell-2p-2n}
    \end{align*}
    which after truncating to $N+1$ terms and performing a Galerkin step with $P_\eta$ yields
    \begin{align}\label{u3expansion}
            \sum_{\ell=0}^N \sum_{m=0}^\ell \sum_{j=0}^m \sum_{p=0}^{j\wedge (m-j)} \sum_{n=0}^{(\ell-m)\wedge (m-2p)} u_{\ell-m} u_j u_{m-j} C(j,m-j,p)C(\ell-m,m-2p,n) \mathbbm{1}_{\{\eta = \ell-2p-2n\}}.
    \end{align}
    This amounts to a number $\tilde N$ of summands shown in Table \ref{Ntable} for low $N$, of which $\tilde N/(N+1)$ are nonzero. For equations with a $u^3$ term, this has a severe impact on the runtime of iPCE for large $N$ (see Remark \ref{multidimensionalremark}).
    \begin{table}[h]
        \centering
        \begin{tabular}{c|ccccccccc}
        \toprule
            $N$ & 0 & 1 & 2 & 3 & 4 & 5 & 6 & 7 & 8 \\
            \midrule
            $\tilde N$ & 1 & 8 & 39 & 124 & 335 & 762 & 1589 & 3016 & 5418 \\ 
            \bottomrule
        \end{tabular}
        \caption{Number $\tilde N$ of summands in \eqref{u3expansion} for low $N$}
        \label{Ntable}
    \end{table}
\end{remark}
 \begin{remark}\label{multidimensionalremark}
    If the random input $\boldsymbol{\xi}$ is multivariate with dimension $P\in\N$, the underlying base of orthogonal polynomial has, for up to degree $N$, $\hat N+1:= \frac{(N+P)!}{N!P!}$ elements \cite[Eq. (5.3)]{xiukarniadakis} (where $\hat N$ denotes the number of expansion terms corresponding to non-constant polynomials). This means that the system \eqref{intrusivepcesystem} is $\hat N+1$-times bigger than the deterministic system. We pointed out in Remark \ref{implementationremark} that especially equations involving a $u^n$ term for $n\geq 3$ can start to pose problems due to the rapidly increasing number of necessary summands in the gPCE. Figure \ref{runtimesN} shows how the runtimes increase for iPCE depending on $N$ for the equations \eqref{differentialequation2} investigated in Section \ref{section:numericalexperiments}. Especially the equation involving a $u^3$ term shows a rapid increase in computational effort as $N$ grows. \\
    On the other hand, suppose that a one-dimensional quadrature formula (such as Gauss-Legendre) has an error bound for $f\in \mathcal{C}^r([-1,1])$ for $q\in \mathbb{N}$ quadrature nodes
    \begin{align*}
        |E_q^1 f| = \mathcal{O}(q^{-r}).
    \end{align*}
    Then, it is known that in $d$ dimensions, the full grid quadrature error $|\tilde E_q^d f|$ and the error for Smolyak sparse quadrature $|E_q^d f|$
    for $f\in H^r([-1,1]^d)$ is given by
    \begin{align*}
       |\tilde E_q^d f| = \mathcal{O}(q^{-r/d}),\quad  |E_q^d f| = \mathcal{O}(\log(q)^{(d-1)(r+1)} q^{-r})
    \end{align*}
    (see e.g. \cite[p. 39]{smolyakmaster}).
    This means that in higher dimensions, provided that the integrand is sufficiently smooth, the number of necessary function evaluations can be drastically reduced by using a sparse grid, and the curse of dimensionality is much less severe than in the iPCE case as outlined above.
\end{remark}
  We will henceforth use the following notation: For $\eta=0,\dots,N$ and $n=0,\dots,M$, it is 
  \begin{align*}
      \boldsymbol{u}_\eta^n := (u_\eta(x_0,t_n),u_\eta(x_1,t_n),\dots,u_\eta(x_p,t_n))^\top.
  \end{align*}
 Also, $\Delta \boldsymbol{u}_\eta^n = (\Delta u_\eta(x_i,t_n))^\top_{i=0,\dots,p}$ refers to $\Delta u_\eta$ evaluated at the grid points.
 \subsubsection{Explicit Euler}\label{section:intrusive_expliciteuler}
 For a single function $u_\eta$ in the PCE, 
 the explicit Euler scheme for solving \eqref{differentialequation} is given by
 \begin{align*}
     \boldsymbol{u}_\eta^{n+1} = \boldsymbol{u}_\eta^n + k\cdot \left[ DA_p\boldsymbol{u}_\eta^n + \sum_{i,j,k=0}^N \boldsymbol{K}_{ijk\eta}\boldsymbol{u}^n_i\odot \boldsymbol{u}^n_j \odot \boldsymbol{u}^n_k \right],\quad n=0,\dots,M-1,
 \end{align*}
 where by $\odot$ we denote a componentwise product. Combining all coefficient functions $\boldsymbol{u}_0,\dots,\boldsymbol{u}_N$, we have
 \begin{algorithm}[iPCE Explicit Euler scheme]\label{algorithm:iPCE_EE}
 \begin{align}
     \boldsymbol{U}^{n+1} &= \boldsymbol{U}^n + k\cdot \left[ (I_{N+1}\otimes DA_p) \boldsymbol{U}^n + \boldsymbol{F}(\boldsymbol{U}^n) \right], \\
     \label{nonlinearfunctionpce}
     \boldsymbol{F}(\boldsymbol{U}^n) &= \left[ \sum_{i,j,k=0}^N \boldsymbol{K}_{ijk\eta} \boldsymbol{u}_i^n\odot \boldsymbol{u}_j^n \odot \boldsymbol{u}_k^n\right]_{\eta=0,\dots,N},\quad n=0,\dots,M-1.
 \end{align}
 \end{algorithm}
 \subsubsection{ETD-RDP and ETD-RDP-IF}\label{section:intrusive_ETDRDPIF}
 Applying the scheme \eqref{etdrdpscheme} to the PCE system \eqref{aftergalerkin} with a finite difference discretization and making use of the fact that $\left(I_{N+1}\otimes (I_p+\frac{k}{3}DA_p)\right)^{-1} = I_{N+1}\otimes (I_p+\frac{k}{3}DA_p)^{-1}$, the resulting scheme is
 \begin{algorithm}[iPCE ETD-RDP scheme]\label{algorithm:iPCE_ETDRDP}
 \begin{align}\label{etdrdp1}
     \boldsymbol{U}^{n+1} &= \left( I_{N+1}\otimes (I_p+\frac{k}{3}DA_p)^{-1}\right)[9\boldsymbol{U}^n + 2k\boldsymbol{F}(\boldsymbol{U}^n) + k\boldsymbol{F}(\boldsymbol{U}^{n+1}_*)] \\ \nonumber
     &- \left( I_{N+1}\otimes (I_p+\frac{k}{4}DA_p)^{-1}\right)[8\boldsymbol{U}^n + \frac{3k}{2}\boldsymbol{F}(\boldsymbol{U}^n) + \frac{k}{2}\boldsymbol{F}(\boldsymbol{U}^{n+1}_*)],\\ \label{etdrdp3}
     \boldsymbol{U}^{n+1}_* &= \left( I_{N+1}\otimes (I_p+kDA_p)^{-1}\right)[\boldsymbol{U}^n + k\boldsymbol{F}(\boldsymbol{U}^n) ],\quad n=0,\dots,M-1
 \end{align}    
 \end{algorithm}
 where $\boldsymbol{F}(\boldsymbol{U}^n)$ is given by \eqref{nonlinearfunctionpce}.\\
 We now discuss the application of ETD-RDP-IF in the case where the spatial dimension is greater than one. For simplicity of the presentation, we restrict ourselves here to the case of two spatial dimensions. In this case, the finite difference discretization of the Laplacian with periodic boundary conditions is given by $A = A_1 + A_2 := (D\cdot I_p) \otimes A_p + A_p \otimes (D\cdot I_p)$. We use the notation $\boldsymbol{U}^n = (\boldsymbol{u}_0,\dots,\boldsymbol{u}_N)^\top$ from above and denote
 \begin{align*}
     \A := I_{N+1}\otimes A,\quad \A_1 := I_{N+1}\otimes A_1,\quad \A_2 := I_{N+1}\otimes A_2,\quad \boldsymbol{I} := I_{N+1}\otimes I_{p^2} = I_{p^2\cdot (N+1)},
 \end{align*}
 so in particular it is $\A = \A_1 + \A_2$. In order to solve the discretized PCE system
 \begin{align}\label{pcediscretized}
     \frac{\partial \boldsymbol{U}}{\partial t} + \A\boldsymbol{U} = \boldsymbol{F}(\boldsymbol{U}),\quad \boldsymbol{u}_0^0 = \boldsymbol{u}_{\mathrm{init}},
 \end{align}
 with $\boldsymbol{F}$ as above, we apply a dimensional splitting technique. We introduce the new time-dependent function $\boldsymbol{V} =  \e^{\A_1t}\boldsymbol{U}$. The term $I_{N+1} \otimes \e^{\A_1t}$ is called the integrating factor for PCE. We will also use that $\A$ and $\A_1$ commute (see \cite[Lemma 1]{Kleefeld2020}), and therefore also $\A$ and $\e^{\A_1 t}$ commute. For the dimensional splitting, we define $\boldsymbol{V}(t):= \e^{\A_1 t} \boldsymbol{U}(t)$. Carrying out the steps as in \cite[pp. 3]{Kleefeld2020} for the PCE case, the time derivative of $\boldsymbol{V}$ is
 \begin{align*}
     \frac{\partial \boldsymbol{V}}{\partial t} = \e^{\A_1 t} \frac{\partial \boldsymbol{U}}{\partial t} + \A_1\e^{\A_1t}\boldsymbol{U}
 \end{align*}
 and inserting \eqref{pcediscretized}, we obtain
 \begin{align*}
     \frac{\partial \boldsymbol{V}}{\partial t} &= \e^{\A_1t}(\boldsymbol{F}(\boldsymbol{U})-\A\boldsymbol{U}) + \A_1 \e^{\A_1t}\boldsymbol{U} = \e^{\A_1t}\boldsymbol{F}(\boldsymbol{U}) - \e^{\A_1 t}\A\boldsymbol{U} + \A_1\e^{\A_1 t}\boldsymbol{U} \\
     &= \e^{\A_1 t}\boldsymbol{F}(\boldsymbol{U})-\A \e^{\A_1 t}\boldsymbol{U} + \A_1 \e^{\A_1 t}\boldsymbol{U}  
     = \e^{\A_1 t}\boldsymbol{F}(\boldsymbol{U}) - \A_2\e^{\A_1 t}\boldsymbol{U} = \e^{\A_1 t}\boldsymbol{F}(\e^{-\A_1 t}\boldsymbol{V}) - \A_2 \boldsymbol{V},
 \end{align*}
 so the system
 \begin{align*}
     \frac{\partial \boldsymbol{V}}{\partial t} + \A_2\boldsymbol{V} = \boldsymbol{G}(\boldsymbol{V}),\quad \boldsymbol{v}_0^0 = \boldsymbol{u}_{\mathrm{init}}
 \end{align*}
 with $\boldsymbol{G}(\boldsymbol{V}) := \e^{\A_1 t}\boldsymbol{F}(\e^{-\A_1 t}\boldsymbol{V})$ must be solved. The next step is to apply the ETD-RDP scheme given in \eqref{etdrdp3} to $\boldsymbol{V}$. The derivation in the PCE case works analogously as described in \cite[p. 5--7]{Kleefeld2020}. In addition to the ETD-RDP scheme \eqref{etdrdp3}, the unwinding of the integrating factor is needed. With the substitutions $\boldsymbol{V}^n = \e^{\A_2 nk}\boldsymbol{U}^n$, $\boldsymbol{V}^{n+1} = \e^{\A_2 nk}\e^{\A_2 k} \boldsymbol{U}^{n+1}$, $\boldsymbol{G}(\boldsymbol{V}^n) = \e^{\A_2 nk} \boldsymbol{F}(\boldsymbol{U}^n)$ and $\boldsymbol{G}(\boldsymbol{V}^{n+1}) = \e^{\A_2 nk} \e^{\A_2 k} \boldsymbol{F}(\boldsymbol{U}^{n+1})$, \eqref{etdrdp3} becomes the full ETD-RDP-IF scheme for PCE:
 \begin{algorithm}[iPCE ETD-RDP-IF scheme]\label{algorithm:iPCE_ETDRDPIF}
 \begin{align*}
     \boldsymbol{U}^{n+1} &= \left( \boldsymbol{I}+\frac{k}{3}D\A_2 \right)^{-1}\left[ \left\{ 9\left(\boldsymbol{I}+\frac{k}{3}D\A_1\right)^{-1} - 8\left( \boldsymbol{I}+\frac{k}{4}D\A_1\right)^{-1} \right\}\right. \\
     &\cdot \left. \{ 9\boldsymbol{U}^n + 2k\boldsymbol{F}(\boldsymbol{U}^n)\} + k\boldsymbol{F}(\boldsymbol{U}_*^{n+1})\right] \\
     &- \left( \boldsymbol{I}+\frac{k}{4}D\A_2 \right)^{-1}\left[ \left\{ 9\left(\boldsymbol{I}+\frac{k}{3}D\A_1\right)^{-1} - 8\left( \boldsymbol{I}+\frac{k}{4}D\A_1\right)^{-1} \right\}\right. \\
     &\cdot \left. \{ 8\boldsymbol{U}^n + \frac{3k}{2}\boldsymbol{F}(\boldsymbol{U}^n)\} + \frac{k}{2}\boldsymbol{F}(\boldsymbol{U}_*^{n+1})\right],\\
     \boldsymbol{U}_*^{n+1} &= (\boldsymbol{I}+kD\A_2)^{-1}(\boldsymbol{I}+kD\A_1)^{-1}(\boldsymbol{U}^n + k\boldsymbol{F}(\boldsymbol{U}^n)),\quad n=0,\dots,M-1.
 \end{align*}
 \end{algorithm}
 Note that this scheme looks identical to the one presented in \cite[Eq. (19)]{Kleefeld2020}, but we have a different notation for $\A_1$, $\A_2$ and the function $\boldsymbol{F}$.
 
 \subsubsection{ETDRK4}\label{section:intrusive_etdrk4}
 For $\eta=0,\dots, N$, we denote $\widehat{\boldsymbol{u}_\eta} = \mathcal{F}(\boldsymbol{u}_\eta)$ and for the stacked discretized PCE functions, we write
 \begin{align*}
     \boldsymbol{U}^n = \begin{pmatrix}
         \boldsymbol{u}_0^n \\
         \vdots \\
         \boldsymbol{u}_N^n
     \end{pmatrix},\quad \widehat{\boldsymbol{U}^n} = \mathcal{F}_{\mathrm{comp}}(\boldsymbol{U}^n) := 
     \begin{pmatrix}
         \mathcal{F}(\boldsymbol{u}_0^n) \\
         \vdots \\
         \mathcal{F}(\boldsymbol{u}_N^n)
     \end{pmatrix}
     =
     \begin{pmatrix}
         \widehat{\boldsymbol{u}_0^n} \\
         \vdots \\
         \widehat{\boldsymbol{u}_N^n}
     \end{pmatrix},\quad 
     \mathcal{F}_{\mathrm{comp}}^{-1}(\widehat{\boldsymbol{U}^n}) := 
     \begin{pmatrix}
         \mathcal{F}^{-1}(\widehat{\boldsymbol{u}_0^n}) \\
         \vdots \\
         \mathcal{F}^{-1}(\widehat{\boldsymbol{u}_N^n})
     \end{pmatrix}
 \end{align*}
 and we stress the fact that the discrete Fourier transform is taken separately for each base function. For the vector $\boldsymbol{U}^n$ of stacked discretized PCE base functions, $\boldsymbol{F}(\boldsymbol{U}^n)$ is understood in the sense of \eqref{nonlinearfunctionpce}. We also use the matrix $L$ from Section \ref{section:etdrk4} and denote $\boldsymbol{L}:= I_{N+1}\otimes L$. This matrix is still diagonal and therefore very easy to invert. As before in \eqref{etdrk4_deterministic}, we introduce the shorthand $\hat{\boldsymbol{F}}(\boldsymbol{\hat U}) := \mathcal{F}_{\mathrm{comp}}(\boldsymbol{F}(\mathcal{F}_{\mathrm{comp}}^{-1}(\hat{\boldsymbol{U}})))$. The ETDRK4 scheme for the intrusive PCE system is then given by
 \begin{algorithm}[iPCE ETDRK4 scheme]\label{algorithm:iPCE_ETDRK4}
 \begin{align*}
     \widehat{\boldsymbol{U}^0} &= \mathcal{F}_{\mathrm{comp}}(\boldsymbol{U}^0), \\
    \widehat{\boldsymbol{U}^{k+1}} &= \exp\left(\frac{h}{2}\boldsymbol{L}\right) \widehat{\boldsymbol{U}^k} \\
    &\ + (h^{-2}\boldsymbol{L}^{-3})\left\{ \left[ -4\boldsymbol{I}-h\boldsymbol{L} + \e^{h\boldsymbol{L}}(4\boldsymbol{I}-3h\boldsymbol{L}+h^2\boldsymbol{L}^2)\right]\hat{\boldsymbol{F}}(\widehat{\boldsymbol{U}^k}) \right. \\
    &\ + 2[2\boldsymbol{I}+h\boldsymbol{L} + \e^{h\boldsymbol{L}}(-2\boldsymbol{I}+h\boldsymbol{L})]\left( \hat{\boldsymbol{F}}(\hat{\boldsymbol{a}}^k) + \hat{\boldsymbol{F}}(\hat{\boldsymbol{b}}^k)\right) \\
    &\ \left. + \left[ -4\boldsymbol{I}-3h\boldsymbol{L} -h^2\boldsymbol{L}^2+\e^{h\boldsymbol{L}}(4\boldsymbol{I}-h\boldsymbol{L})\right] \hat{\boldsymbol{F}}(\hat{\boldsymbol{c}}^k) \right\},\\
    {\boldsymbol{U}}^M &= \mathcal{F}_{\mathrm{comp}}^{-1}(\widehat {\boldsymbol{U}^M}),
 \end{align*}
 \end{algorithm}
 where $\boldsymbol{I} : = I_{(N+1)p^2}$ and
\begin{align*}
    \hat{\boldsymbol{a}}^k &= \e^{\boldsymbol{L}h/2}\widehat{\boldsymbol{U}^k} + \boldsymbol{L}^{-1}(\e^{\boldsymbol{L}h/2}-\boldsymbol{I})\hat{\boldsymbol{F}}(\widehat{\boldsymbol{U}^k}),\\
    \hat{\boldsymbol{b}}^k &= \e^{\boldsymbol{L}h/2}\widehat{\boldsymbol{U}^k} + \boldsymbol{L}^{-1}(\e^{\boldsymbol{L}h/2}-\boldsymbol{I})\hat{\boldsymbol{F}}(\hat{\boldsymbol{a}}^k), \\
    \hat{\boldsymbol{c}}^k &= \e^{\boldsymbol{L}h/2}\hat{\boldsymbol{a}}^k + \boldsymbol{L}^{-1}(\e^{\boldsymbol{L}h/2}-\boldsymbol{I})(2\hat{\boldsymbol{F}}(\hat{\boldsymbol{b}}^k)-\hat{\boldsymbol{F}}(\widehat{\boldsymbol{U}^k})).
\end{align*}
In Section \ref{section:etdrk4}, we recalled a contour integral method in order to avoid cancellation errors arising in the evaluation of the expressions $a_n$, $b_n$ and $c_n$ in equations \eqref{deterministic_an} to \eqref{deterministic_cn}. In this scheme, we apply this technique to each base function, i.e. the mean \eqref{complexrootsmean} is computed separately for each PCE base function. Likewise, anti-aliasing is applied to each base function.
\begin{remark}
    In order to deal with random PDEs of the form 
    \begin{align*}
        \frac{\partial u(x,t,\omega)}{\partial t} = D(\omega) \Delta u(x,t,\omega) + F(u(x,t,\omega)),
    \end{align*}
    i.e. equations where $D$ is random and the function $F$ is not random, the ETDRK4 intrusive PCE scheme described above is not very suitable, because in this case the discretized Laplacian $D\Delta$ is not a tridiagonal matrix anymore (the shape of the matrix depends on the distribution of the random variable $D$), and the Fourier transform $\boldsymbol{L}$ of the Laplacian is not diagonal anymore, which can severely impact the performance of the scheme.
\end{remark}
\begin{figure}
\newcommand{\X}{2.5}
\newcommand{\Yone}{1.0}
\newcommand{\Ytwo}{0.5}
\newcommand{\intra}{runtimearray_EE.txt}
\newcommand{\intrb}{runtimearray_ETDRDP.txt}
\newcommand{\intrc}{runtimearray_ETDRK4.txt}
\newcommand{\legendposition}{north west}
\centering

\begin{minipage}{\textwidth} 
\centering 
\captionsetup{justification=centering} 
\caption*{\underline{Runtimes for iPCE for different $N$ compared to $N=0$}}
\end{minipage}

 \begin{tabular}{@{}c@{}}
\begin{tikzpicture}
    \begin{axis}
    [ cycle list name=goodcolors,
    title={\textbf{EE}},
      height=\Heightt, width=\Widtht,
      xlabel={$N$}, ylabel={relative runtime},ymode = log,
      xtick={0,3,6,9},
      legend entries={Lin. $D=0$, Lin. $D=1$},
      legend pos=north west, ]
      \addlegendimage{color=plotcolor1,solid, line width=2pt}
    \addlegendimage{color=plotcolor2,dashed, line width=2pt}
    \foreach \i in {1,2,3,4,5,6} {
    \addplot table [col sep=comma, x index=0, y index=\i] {\intra};
    }
    \end{axis}
\end{tikzpicture}
\begin{tikzpicture}
    \begin{axis}
    [ cycle list name=goodcolors,
    title={\textbf{ETD-RDP}},
      height=\Heightt, width=\Widtht,
      xlabel={$N$},ymode = log,
      xtick={0,3,6,9},
      legend entries={Quad. $D=0$,Quad. $D=1$},
      legend pos=north west, ]
        \addlegendimage{plotcolor3,dotted, line width=2pt}
        \addlegendimage{plotcolor4,densely dashed, line width=2pt}
    \foreach \i in {1,2,3,4,5,6} {
    \addplot table [col sep=comma, x index=0, y index=\i] {\intrb};
    }
    \end{axis}
\end{tikzpicture}
\begin{tikzpicture}
    \begin{axis}
    [ cycle list name=goodcolors,
    title={\textbf{ETDRK4}},
      height=\Heightt, width=\Widtht,
      xlabel={$N$},ymode = log,
      xtick={0,3,6,9},
      legend entries={Cub. $D=0$,Cub. $D=1$},
      legend pos=north west, ]
      \addlegendimage{plotcolor5,dashed, dash pattern=on 6pt off 2pt, line width=2pt}
        \addlegendimage{plotcolor6,solid, line width=1pt}
    \foreach \i in {1,2,3,4,5,6} {
    \addplot table [col sep=comma, x index=0, y index=\i] {\intrc};
    }
    \end{axis}
\end{tikzpicture}
\end{tabular}
\caption{Relative runtimes $R_N$ for $N=0,\dots,9$, where $R_0:=1$. Each plot shows six curves for six different equations investigated in Section \ref{section:numericalexperiments}: Equation \eqref{differentialequation2} with linear, quadratic or cubic $F$, with $D=0$ or $D=1$. Each data point is an average over the runtimes of ten identical iPCE simulations with $T=0.1$ and $M=10$ time steps.} \label{runtimesN}
\end{figure}
\section{Numerical experiments.}\label{section:numericalexperiments}
We now apply non-intrusive PCE as explained in Section \ref{section:nonintrusivepce} and intrusive PCE as explained in Section \ref{section:intrusivepce} to the equation \eqref{differentialequation} which we repeat here:
\begin{align}\label{differentialequation2}
\begin{cases}
    \frac{\partial u(x,t,\omega)}{\partial t} &= D\Delta u(x,t,\omega) + F(\omega,u(x,t,\omega)),\\
    u(x,0,\omega) &= u_\mathrm{init}(x).
\end{cases}
\end{align}
For a random constant $K\sim \mathcal{U}[1,2]$, we will use $F(\omega,u(x,t,\omega)) = K(\omega) u(x,t,\omega)$ in Section \ref{section:linearequation}, $F(\omega,u(x,t,\omega)) = K(\omega)u(x,t,\omega)^2$ in Section \ref{section:quadraticequation} and $F(\omega,u(x,t,\omega)) = K(\omega)u(x,t,\omega)^3$ in Section \ref{section:cubicequation}. In all three cases, we test the two cases $D=1$ and $D=0$, i.e. a random ODE and a PDE with a diffusion term. In all three cases, we pick $u_{\mathrm{init}}(x) = \cos(\pi x)$.  Finally, we will discuss a more complex problem, the Gray-Scott system, in Section \ref{section:grayscott}. Throughout all simulations, we work on the spatial domain $(-1,1)^d$, and spatial resolution $p=128$. \\
In all error plots, we show the relative $L^2$ error which, for the exact mean or a reference solution mean $\mathbb{E}[u]$ and the approximated solution $\mathbb{E}[u]_{\mathrm{approx}}$, is given by $e_{\mathrm{rel}}(t) = \|\mathbb{E}[u](\cdot,t)-\mathbb{E}[u]_{\mathrm{approx}}(\cdot,t)\|_{L^2}/\|\mathbb{E}[u](\cdot,t)\|_{L^2}$. For the random equation with linear term in Section \ref{section:linearequation}, we know the exact solution, for the other equations in Sections \ref{section:quadraticequation}, \ref{section:cubicequation} and \ref{section:grayscott} we use a reference solution obtained from using non-intrusive PCE with Gauß quadrature with a high number of quadrature points (see Table \ref{runtimetable}). We also note that for the error of niPCE with naive, randomly sampled MC, we ran ten simulations producing ten error graphs over which we took the mean. 
\subsection{One-dimensional random equation with linear term}\label{section:linearequation}
The equation 
\begin{align}\label{linearequation}
    \frac{\partial u(x,t,\omega)}{\partial t} = D\Delta u(x,t,\omega) - K(\omega)u(x,t,\omega),\quad u(x,0,\omega) = \cos(\pi x),\quad x\in (-1,1),\ t\in (0,T]
\end{align}
has the exact solution $u(x,t,\omega) = \exp(-(K(\omega)+D\pi^2)t)\cos(\pi x)$. If $K\sim \mathcal{U}[a,b]$, the expected value is given by
\begin{align} \nonumber
        \mathbb{E}[u(x,t,\cdot)] &= \frac{1}{b-a}\int_a^b \exp(-(\xi+D \pi^2)t)\cos(\pi x)\; \de \xi \\ \label{1D_exactmean}
        &= \frac{1}{b-a}\frac{\cos(\pi x)}{t}  \exp(-(D\pi^2+ a)t)-\exp(-(D\pi^2 + b)t),
\end{align}
and in our tests we choose again $a=1$ and $b=2$.
The variance is 
\begin{align} \nonumber
    \mathrm{Var}[u(x,t,\cdot)] &= \mathbb{E}[u(x,t,\cdot)^2] - \mathbb{E}[u(x,t,\cdot)]^2 \\ \nonumber
    &= \frac{1}{b-a}\int_a^b \exp(-2(\xi+D\pi^2)t)\cos^2(\pi x)\; \de \xi - \mathbb{E}[u(x,t,\cdot)]^2 \\ \nonumber
    &= \frac{1}{b-a}\frac{\cos^2(\pi x)}{2t}\left[ \exp(-2(D\pi^2+ a)t)-\exp(-2(D\pi^2 + b)t)\right] \\ \label{exactvariance}
    &- \left[ \frac{1}{b-a}\frac{\cos(\pi x)}{t}  \exp(-(D\pi^2+ a)t)-\exp(-(D\pi^2 + b)t) \right]^2. 
\end{align}
 In Figure \ref{linear_D=0}, we show the time-dependent error plots for solving \eqref{linearequation} with $D=0$, comparing to the exact solution. For iPCE, it can be seen that both for EE and for ETD-RDP, increasing the number of non-constant polynomials beyond $N=1$ makes no difference. Due to the greater accuracy of ETDRK4, the error keeps decreasing until $N=4$, and the spectral convergence of iPCE is clearly seen. In the plots for niPCE, it can be seen that the error strongly depends on the method of numerical quadrature to compute \eqref{pcenumint}. By choice of the quadrature points and weights alone, niPCE with GQ achieves better accuracy than iPCE in all scenarios. While, therefore, niPCE is a clear winner in this case, the situation is different for $D=1$, where it can be seen in Figure \ref{linear_D=1} that for ETD-RDP and ETDRK4, a lower error can be achieved with iPCE for comparable runtimes (see Table \ref{runtimetable}).\\
In Figure \ref{linear_D=0_variance}, the error plots for the variance are shown, computed according to formula \eqref{meanandvariance}. In the case of iPCE, the iPCE scheme is run and $\sum_{i=1}^N|u_i(x,T)|^2$ is computed from the result of the scheme. In the case of niPCE, formula \eqref{pcenumint} is applied with $q=10$. The EE, ETD-RDP and ETDRK4 schemes again form a clear hierarchy from least to most accurate. It can also be seen that the error for ETDRK4 iPCE, $N=5$, is very similar to that of ETDRK4 niPCE, indicating that the remaining error might not be due to the choice of the numerical scheme ETDRK4 but rather due to approximation steps regarding the gPCE, such as in \eqref{u3expansion}. Also, MC and QMC yield high errors, which is due to the fact that the norm of the variance is very low compared to the mean, which makes the sampling error in \eqref{pcenumint} much more significant. For instance, for the canonical normed Legendre polynomial $P_1$, for a random sequence $\{x^{\mathrm{MC}}_j\}_{j=1,\dots,100}$, a Sobol sequence $\{x^{\mathrm{QMC}}_j\}_{j=1,\dots,100}$ and Gauss-Legendre quadrature nodes $\{x^{\mathrm{GQ}}_j\}_{j=1,\dots,100}$ it is (with the corresponding weights as described in \eqref{pcenumint})
\begin{align}\label{samplingerrorexample}
    \sum_{j=1}^{100} w_j^{\mathrm{MC}} P_1(x^{\mathrm{MC}}_j) \approx -0.00898,\;\; \sum_{j=1}^{100} w_j^{\mathrm{QMC}} P_1(x^{\mathrm{QMC}}_j) \approx -0.00622, \;\; \sum_{j=1}^{100} w_j^{\mathrm{GQ}} P_1(x^{\mathrm{GQ}}_j) \approx 2.02\cdot 10^{-17}
\end{align}
(the first number is an average over 100 Monte Carlo samples) due to MC and QMC being non-symmetric rules, while Gauss-Legendre quadrature nodes are symmetric about the origin. In particular for small $t$, due to the deterministic initial condition the variance is very close to zero, which, due to sampling errors as demonstrated in \eqref{samplingerrorexample}, leads to a very large relative error.
\begin{figure}
\newcommand{\X}{2.5}
\newcommand{\Yone}{1.0}
\newcommand{\Ytwo}{0.5}
\newcommand{\intra}{errorarray_FD_EE_mean_system=6_D=0.00000.txt}
\newcommand{\intrb}{errorarray_FD_ETDRDP_mean_system=6_D=0.00000.txt}
\newcommand{\intrc}{errorarray_Spectral_mean_system=6_D=0.00000.txt}
\newcommand{\nintra}{errorarray_mean_Nonintrusive_FD_EE_system=6_D=0.00000.txt}
\newcommand{\nintrb}{errorarray_mean_Nonintrusive_FD_ETDRDP_system=6_D=0.00000.txt}
\newcommand{\nintrc}{errorarray_mean_Nonintrusive_Spectral_system=6_D=0.00000.txt}
\newcommand{\legendposition}{north west}
\centering

\begin{minipage}{\textwidth} 
\centering 
\captionsetup{justification=centering} 
\caption*{\underline{Equation with linear term, $D=0$, mean}}
\end{minipage}
 \begin{tabular}{@{}c@{}}
\begin{tikzpicture}
    \begin{axis}
    [ cycle list name=goodcolors,
    title={Intrusive PCE, \textbf{EE}},
      height=\Height, width=\Width,
      xlabel={$t$}, ylabel={relative $L^2$ error},ymode = log,
      legend entries={N=1, N=2, N=3, N=4, N=5},
      legend pos=south east, ]
    \foreach \i in {1,2,3,4,5} {
    \addplot table [col sep=comma, x index=0, y index=\i] {\intra}; }
    \end{axis}
\end{tikzpicture}
\begin{tikzpicture}
    \begin{axis}
    [ cycle list name=goodcolors,
    title={Non-Intrusive PCE, \textbf{EE}},
      height=\Height, width=\Width,
      xlabel={$t$}, ylabel={relative $L^2$ error},ymode = log,
      legend entries={MC, QMC, GQ},
      legend pos=south east, ]
    \foreach \i in {1,2,3} {
    \addplot table [col sep=comma, x index=0, y index=\i] {\nintra}; }
    \end{axis}
\end{tikzpicture}
\\
\begin{tikzpicture}
    \begin{axis}
    [ cycle list name=goodcolors,
    title={Intrusive PCE, \textbf{ETD-RDP}},
      height=\Height, width=\Width,
      xlabel={$t$}, ylabel={relative $L^2$ error},ymode = log,
      legend pos=\legendposition, ]
    \foreach \i in {1,2,3,4,5} {
    \addplot table [col sep=comma, x index=0, y index=\i] {\intrb}; }
    \end{axis}
\end{tikzpicture}
\begin{tikzpicture}
    \begin{axis}
    [ cycle list name=goodcolors,
    title={Non-Intrusive PCE, \textbf{ETD-RDP}},
      height=\Height, width=\Width,
      xlabel={$t$}, ylabel={relative $L^2$ error},ymode = log,
      legend pos=\legendposition, ]
    \foreach \i in {1,2,3} {
    \addplot table [col sep=comma, x index=0, y index=\i] {\nintrb}; }
    \end{axis}
\end{tikzpicture}
\\
\begin{tikzpicture}
    \begin{axis}
    [ cycle list name=goodcolors,
    title={Intrusive PCE, \textbf{ETDRK4}},
      height=\Height, width=\Width,
      xlabel={$t$}, ylabel={relative $L^2$ error},ymode = log,
      legend pos=\legendposition, ]
    \foreach \i in {1,2,3,4,5} {
    \addplot table [col sep=comma, x index=0, y index=\i] {\intrc}; }
    \end{axis}
\end{tikzpicture}
\begin{tikzpicture}
    \begin{axis}
    [ cycle list name=goodcolors,
    title={Non-Intrusive PCE, \textbf{ETDRK4}},
      height=\Height, width=\Width,
      xlabel={$t$}, ylabel={relative $L^2$ error},ymode = log,
      legend pos=\legendposition, ]
    \foreach \i in {1,2,3} {
    \addplot table [col sep=comma, x index=0, y index=\i] {\nintrc}; }
    \end{axis}
\end{tikzpicture}
\end{tabular}
\caption{Plots showing the relative time-dependent error for iPCE schemes (Algorithms \ref{algorithm:iPCE_EE}, \ref{algorithm:iPCE_ETDRDP}, \ref{algorithm:iPCE_ETDRK4}, left-hand side) and niPCE schemes (right-hand side) for equation \eqref{linearequation} with linear term with diffusion constant $D=0$. The exact mean is given by \eqref{1D_exactmean}.} \label{linear_D=0}
\end{figure}
\begin{figure}
\newcommand{\X}{2.5} \newcommand{\Yone}{1.0}
\newcommand{\Ytwo}{0.5} \newcommand{\intra}{errorarray_FD_EE_mean_system=6_D=1.00000.txt}
\newcommand{\intrb}{errorarray_FD_ETDRDP_mean_system=6_D=1.00000.txt}
\newcommand{\intrc}{errorarray_Spectral_mean_system=6_D=1.00000.txt}
\newcommand{\nintra}{errorarray_mean_Nonintrusive_FD_EE_system=6_D=1.00000.txt}
\newcommand{\nintrb}{errorarray_mean_Nonintrusive_FD_ETDRDP_system=6_D=1.00000.txt}
\newcommand{\nintrc}{errorarray_mean_Nonintrusive_Spectral_system=6_D=1.00000.txt}
\newcommand{\legendposition}{north west}
\centering

\begin{minipage}{\textwidth} 
\centering 
\captionsetup{justification=centering} 
\caption*{\underline{Equation with linear term, $D=1$, mean}}
\end{minipage}
 \begin{tabular}{@{}c@{}}
\begin{tikzpicture}
    \begin{axis}
    [ cycle list name=goodcolors,
    title={Intrusive PCE, \textbf{EE}},
      height=\Height, width=\Width,
      xlabel={$t$}, ylabel={relative $L^2$ error},ymode = log,
      legend entries={N=1, N=2, N=3, N=4, N=5},
      legend pos=south east, ]
    \foreach \i in {1,2,3,4,5} {
    \addplot table [col sep=comma, x index=0, y index=\i] {\intra}; }
    \end{axis}
\end{tikzpicture}
\begin{tikzpicture}
    \begin{axis}
    [ cycle list name=goodcolors,
    title={Non-Intrusive PCE, \textbf{EE}},
      height=\Height, width=\Width,
      xlabel={$t$}, ylabel={relative $L^2$ error},ymode = log,
      legend entries={MC, QMC, GQ},
      legend pos=south east, ]
    \foreach \i in {1,2,3} {
    \addplot table [col sep=comma, x index=0, y index=\i] {\nintra}; }
    \end{axis}
\end{tikzpicture}
\\
\begin{tikzpicture}
    \begin{axis}
    [ cycle list name=goodcolors,
    title={Intrusive PCE, \textbf{ETD-RDP}},
      height=\Height, width=\Width,
      xlabel={$t$}, ylabel={relative $L^2$ error},ymode = log,
      legend pos=\legendposition, ]
    \foreach \i in {1,2,3,4,5} {
    \addplot table [col sep=comma, x index=0, y index=\i] {\intrb}; }
    \end{axis}
\end{tikzpicture}
\begin{tikzpicture}
    \begin{axis}
    [ cycle list name=goodcolors,
    title={Non-Intrusive PCE, \textbf{ETD-RDP}},
      height=\Height, width=\Width,
      xlabel={$t$}, ylabel={relative $L^2$ error},ymode = log,
      legend pos=\legendposition, ]
    \foreach \i in {1,2,3} {
    \addplot table [col sep=comma, x index=0, y index=\i] {\nintrb}; }
    \end{axis}
\end{tikzpicture}
\\
\begin{tikzpicture}
    \begin{axis}
    [ cycle list name=goodcolors,
    title={Intrusive PCE, \textbf{ETDRK4}},
      height=\Height, width=\Width,
      xlabel={$t$}, ylabel={relative $L^2$ error},ymode = log,
      legend pos=\legendposition, ]
    \foreach \i in {1,2,3,4,5} {
    \addplot table [col sep=comma, x index=0, y index=\i] {\intrc}; }
    \end{axis}
\end{tikzpicture}
\begin{tikzpicture}
    \begin{axis}
    [ cycle list name=goodcolors,
    title={Non-Intrusive PCE, \textbf{ETDRK4}},
      height=\Height, width=\Width,
      xlabel={$t$}, ylabel={relative $L^2$ error},ymode = log,
      legend pos=\legendposition, ]
    \foreach \i in {1,2,3} {
    \addplot table [col sep=comma, x index=0, y index=\i] {\nintrc}; }
    \end{axis}
\end{tikzpicture}
\end{tabular}
\caption{Plots showing the relative time-dependent error for iPCE schemes (Algorithms \ref{algorithm:iPCE_EE}, \ref{algorithm:iPCE_ETDRDP}, \ref{algorithm:iPCE_ETDRK4}, left-hand side) and niPCE schemes (right-hand side) for equation \eqref{linearequation} with linear term with diffusion constant $D=1$. The exact mean is given by \eqref{1D_exactmean}.}\label{linear_D=1}
\end{figure}
\begin{figure}
\newcommand{\X}{2.5}
\newcommand{\Yone}{1.0}
\newcommand{\Ytwo}{0.5}
\newcommand{\intra}{errorarray_FD_EE_variance_system=6_D=0.00000.txt}
\newcommand{\intrb}{errorarray_FD_ETDRDP_variance_system=6_D=0.00000.txt}
\newcommand{\intrc}{errorarray_Spectral_variance_system=6_D=0.00000.txt}
\newcommand{\nintra}{errorarray_variance_Nonintrusive_FD_EE_system=6_D=0.00000.txt}
\newcommand{\nintrb}{errorarray_variance_Nonintrusive_FD_ETDRDP_system=6_D=0.00000.txt}
\newcommand{\nintrc}{errorarray_variance_Nonintrusive_Spectral_system=6_D=0.00000.txt}
\newcommand{\legendposition}{north west}
\centering

\begin{minipage}{\textwidth} 
\centering 
\captionsetup{justification=centering} 
\caption*{\underline{Equation with linear term, $D=0$, variance}}
\end{minipage}

 \begin{tabular}{@{}c@{}}
\begin{tikzpicture}
    \begin{axis}
    [ cycle list name=goodcolors,
    title={Intrusive PCE, \textbf{EE}},
      height=\Height, width=\Width,
      xlabel={$t$}, ylabel={relative $L^2$ error},ymode = log,
      legend entries={N=1, N=2, N=3, N=4, N=5},
      legend pos=north west, ]
    \foreach \i in {1,2,3,4,5} {
    \addplot table [col sep=comma, x index=0, y index=\i] {\intra}; }
    \end{axis}
\end{tikzpicture}
\begin{tikzpicture}
    \begin{axis}
    [ cycle list name=goodcolors,
    title={Non-Intrusive PCE, \textbf{EE}},
      height=\Height, width=\Width,
      xlabel={$t$}, ylabel={relative $L^2$ error},ymode = log,
      legend entries={MC, QMC, GQ},
      legend pos=north east, ]
    \foreach \i in {1,2,3} {
    \addplot table [col sep=comma, x index=0, y index=\i] {\nintra}; }
    \end{axis}
\end{tikzpicture}
\\
\begin{tikzpicture}
    \begin{axis}
    [ cycle list name=goodcolors,
    title={Intrusive PCE, \textbf{ETD-RDP}},
      height=\Height, width=\Width,
      xlabel={$t$}, ylabel={relative $L^2$ error},ymode = log,
      legend pos=\legendposition, ]
    \foreach \i in {1,2,3,4,5} {
    \addplot table [col sep=comma, x index=0, y index=\i] {\intrb}; }
    \end{axis}
\end{tikzpicture}
\begin{tikzpicture}
    \begin{axis}
    [ cycle list name=goodcolors,
    title={Non-Intrusive PCE, \textbf{ETD-RDP}},
      height=\Height, width=\Width,
      xlabel={$t$}, ylabel={relative $L^2$ error},ymode = log,
      legend pos=\legendposition, ]
    \foreach \i in {1,2,3} {
    \addplot table [col sep=comma, x index=0, y index=\i] {\nintrb}; }
    \end{axis}
\end{tikzpicture}
\\
\begin{tikzpicture}
    \begin{axis}
    [ cycle list name=goodcolors,
    title={Intrusive PCE, \textbf{ETDRK4}},
      height=\Height, width=\Width,
      xlabel={$t$}, ylabel={relative $L^2$ error},ymode = log,
      legend pos=\legendposition, ]
    \foreach \i in {1,2,3,4,5} {
    \addplot table [col sep=comma, x index=0, y index=\i] {\intrc}; }
    \end{axis}
\end{tikzpicture}
\begin{tikzpicture}
    \begin{axis}
    [ cycle list name=goodcolors,
    title={Non-Intrusive PCE, \textbf{ETDRK4}},
      height=\Height, width=\Width,
      xlabel={$t$}, ylabel={relative $L^2$ error},ymode = log,
      legend pos=\legendposition, ]
    \foreach \i in {1,2,3} {
    \addplot table [col sep=comma, x index=0, y index=\i] {\nintrc}; }
    \end{axis}
\end{tikzpicture}
\end{tabular}
\caption{Plots showing the relative time-dependent error for the same setup as in Figure \ref{linear_D=0}, but for the variance computed according to equation \eqref{meanandvariance} with $q=10$ coefficient functions using formula \eqref{pcenumint}. The large errors for MC and QMC are due to the variance being close to zero, so that the sampling error introduced in \eqref{pcenumint} plays a much bigger role (see \eqref{samplingerrorexample}). The exact variance is given by \eqref{exactvariance}.}\label{linear_D=0_variance}
\end{figure}
\subsection{One-dimensional random equation with quadratic term}\label{section:quadraticequation}
We consider the equation with a quadratic term
\begin{align}\label{quadraticequation}
    \frac{\partial u(x,t,\omega)}{\partial t} = D\Delta u(x,t,\omega) - K(\omega)u(x,t,\omega)^2,\quad u(x,0,\omega) = \cos(\pi x),\quad x\in (-1,1),\ t\in (0,T].
\end{align}
Numerical simulations show that for $D=0$, this equation's solution diverges as early as about $t=0.5$ for $K=2$, which is why the errors in Figure \ref{quadratic_D=0} also show a sharp increase towards that point in time. In this case, despite the nonlinearity now being quadratic, iPCE is competitive with niPCE for all three schemes. This is also true for the case $D=1$ for ETD-RDP and for ETDRK4, as can be seen in Figure \ref{quadratic_D=1}.

\begin{figure}
\newcommand{\X}{2.5}
\newcommand{\Yone}{1.0}
\newcommand{\Ytwo}{0.5}
\newcommand{\intra}{errorarray_FD_EE_mean_system=7_D=0.00000.txt}
\newcommand{\intrb}{errorarray_FD_ETDRDP_mean_system=7_D=0.00000.txt}
\newcommand{\intrc}{errorarray_Spectral_mean_system=7_D=0.00000.txt}
\newcommand{\nintra}{errorarray_mean_Nonintrusive_FD_EE_system=7_D=0.00000.txt}
\newcommand{\nintrb}{errorarray_mean_Nonintrusive_FD_ETDRDP_system=7_D=0.00000.txt}
\newcommand{\nintrc}{errorarray_mean_Nonintrusive_Spectral_system=7_D=0.00000.txt}
\newcommand{\legendposition}{north west}
\centering
\begin{minipage}{\textwidth} 
\centering 
\captionsetup{justification=centering} 
\caption*{\underline{Equation with quadratic term, $D=0$, mean}}
\end{minipage}
 \begin{tabular}{@{}c@{}}
\begin{tikzpicture}
    \begin{axis}
    [ cycle list name=goodcolors,
    title={Intrusive PCE, \textbf{EE}},
      height=\Height, width=\Width,
      xlabel={$t$}, ylabel={relative $L^2$ error},ymode = log,
      legend entries={N=1, N=2, N=3, N=4, N=5},
      legend pos=south east, ]
    \foreach \i in {1,2,3,4,5} {
    \addplot table [col sep=comma, x index=0, y index=\i] {\intra}; }
    \end{axis}
\end{tikzpicture}
\begin{tikzpicture}
    \begin{axis}
    [ cycle list name=goodcolors,
    title={Non-Intrusive PCE, \textbf{EE}},
      height=\Height, width=\Width,
      xlabel={$t$}, ylabel={relative $L^2$ error},ymode = log,
      legend entries={MC, QMC, GQ},
      legend pos=\legendposition, ]
    \foreach \i in {1,2,3} {
    \addplot table [col sep=comma, x index=0, y index=\i] {\nintra}; }
    \end{axis}
\end{tikzpicture}
\\
\begin{tikzpicture}
    \begin{axis}
    [ cycle list name=goodcolors,
    title={Intrusive PCE, \textbf{ETD-RDP}},
      height=\Height, width=\Width,
      xlabel={$t$}, ylabel={relative $L^2$ error},ymode = log,
      legend pos=\legendposition, ]
    \foreach \i in {1,2,3,4,5} {
    \addplot table [col sep=comma, x index=0, y index=\i] {\intrb}; }
    \end{axis}
\end{tikzpicture}
\begin{tikzpicture}
    \begin{axis}
    [ cycle list name=goodcolors,
    title={Non-Intrusive PCE, \textbf{ETD-RDP}},
      height=\Height, width=\Width,
      xlabel={$t$}, ylabel={relative $L^2$ error},ymode = log,
      legend pos=\legendposition, ]
    \foreach \i in {1,2,3} {
    \addplot table [col sep=comma, x index=0, y index=\i] {\nintrb}; }
    \end{axis}
\end{tikzpicture}
\\
\begin{tikzpicture}
    \begin{axis}
    [ cycle list name=goodcolors,
    title={Intrusive PCE, \textbf{ETDRK4}},
      height=\Height, width=\Width,
      xlabel={$t$}, ylabel={relative $L^2$ error},ymode = log,
      legend pos=\legendposition, ]
    \foreach \i in {1,2,3,4,5} {
    \addplot table [col sep=comma, x index=0, y index=\i] {\intrc}; }
    \end{axis}
\end{tikzpicture}
\begin{tikzpicture}
    \begin{axis}
    [ cycle list name=goodcolors,
    title={Non-Intrusive PCE, \textbf{ETDRK4}},
      height=\Height, width=\Width,
      xlabel={$t$}, ylabel={relative $L^2$ error},ymode = log,
      legend pos=\legendposition, ]
    \foreach \i in {1,2,3} {
    \addplot table [col sep=comma, x index=0, y index=\i] {\nintrc}; }
    \end{axis}
\end{tikzpicture}
\end{tabular}
\caption{Plots showing the relative time-dependent error for iPCE schemes (Algorithms \ref{algorithm:iPCE_EE}, \ref{algorithm:iPCE_ETDRDP}, \ref{algorithm:iPCE_ETDRK4}, left-hand side) and niPCE schemes (right-hand side) for equation \eqref{quadraticequation} with quadratic term with diffusion constant $D=0$.}\label{quadratic_D=0}
\end{figure}
\begin{figure}
\newcommand{\X}{2.5} \newcommand{\Yone}{1.0}
\newcommand{\Ytwo}{0.5} \newcommand{\intra}{errorarray_FD_EE_mean_system=7_D=1.00000.txt}
\newcommand{\intrb}{errorarray_FD_ETDRDP_mean_system=7_D=1.00000.txt}
\newcommand{\intrc}{errorarray_Spectral_mean_system=7_D=1.00000.txt}
\newcommand{\nintra}{errorarray_mean_Nonintrusive_FD_EE_system=7_D=1.00000.txt}
\newcommand{\nintrb}{errorarray_mean_Nonintrusive_FD_ETDRDP_system=7_D=1.00000.txt}
\newcommand{\nintrc}{errorarray_mean_Nonintrusive_Spectral_system=7_D=1.00000.txt}
\newcommand{\legendposition}{north west}
\centering

\begin{minipage}{\textwidth} 
\centering 
\captionsetup{justification=centering} 
\caption*{\underline{Equation with quadratic term, $D=1$, mean}}
\end{minipage}
 \begin{tabular}{@{}c@{}}
\begin{tikzpicture}
    \begin{axis}
    [ cycle list name=goodcolors,
    title={Intrusive PCE, \textbf{EE}},
      height=\Height, width=\Width,
      xlabel={$t$}, ylabel={relative $L^2$ error},ymode = log,
      legend entries={N=1, N=2, N=3, N=4, N=5},
      legend pos=south east, ]
    \foreach \i in {1,2,3,4,5} {
    \addplot table [col sep=comma, x index=0, y index=\i] {\intra}; }
    \end{axis}
\end{tikzpicture}
\begin{tikzpicture}
    \begin{axis}
    [ cycle list name=goodcolors,
    title={Non-Intrusive PCE, \textbf{EE}},
      height=\Height, width=\Width,
      xlabel={$t$}, ylabel={relative $L^2$ error},ymode = log,
      legend entries={MC, QMC, GQ},
      legend pos=south east, ]
    \foreach \i in {1,2,3} {
    \addplot table [col sep=comma, x index=0, y index=\i] {\nintra}; }
    \end{axis}
\end{tikzpicture}
\\
\begin{tikzpicture}
    \begin{axis}
    [ cycle list name=goodcolors,
    title={Intrusive PCE, \textbf{ETD-RDP}},
      height=\Height, width=\Width,
      xlabel={$t$}, ylabel={relative $L^2$ error},ymode = log,
      legend pos=south east, ]
    \foreach \i in {1,2,3,4,5} {
    \addplot table [col sep=comma, x index=0, y index=\i] {\intrb}; }
    \end{axis}
\end{tikzpicture}
\begin{tikzpicture}
    \begin{axis}
    [ cycle list name=goodcolors,
    title={Non-Intrusive PCE, \textbf{ETD-RDP}},
      height=\Height, width=\Width,
      xlabel={$t$}, ylabel={relative $L^2$ error},ymode = log,
      legend pos=\legendposition, ]
    \foreach \i in {1,2,3} {
    \addplot table [col sep=comma, x index=0, y index=\i] {\nintrb}; }
    \end{axis}
\end{tikzpicture}
\\
\begin{tikzpicture}
    \begin{axis}
    [ cycle list name=goodcolors,
    title={Intrusive PCE, \textbf{ETDRK4}},
      height=\Height, width=\Width,
      xlabel={$t$}, ylabel={relative $L^2$ error},ymode = log,
      legend pos=\legendposition, ]
    \foreach \i in {1,2,3,4,5} {
    \addplot table [col sep=comma, x index=0, y index=\i] {\intrc}; }
    \end{axis}
\end{tikzpicture}
\begin{tikzpicture}
    \begin{axis}
    [ cycle list name=goodcolors,
    title={Non-Intrusive PCE, \textbf{ETDRK4}},
      height=\Height, width=\Width,
      xlabel={$t$}, ylabel={relative $L^2$ error},ymode = log,
      legend pos=\legendposition, ]
    \foreach \i in {1,2,3} {
    \addplot table [col sep=comma, x index=0, y index=\i] {\nintrc}; }
    \end{axis}
\end{tikzpicture}
\end{tabular}
\caption{Plots showing the relative time-dependent error for iPCE schemes (Algorithms \ref{algorithm:iPCE_EE}, \ref{algorithm:iPCE_ETDRDP}, \ref{algorithm:iPCE_ETDRK4}, left-hand side) and niPCE schemes (right-hand side) for equation \eqref{quadraticequation} with quadratic term with diffusion constant $D=1$.}\label{quadratic_D=1}
\end{figure}
\subsection{One-dimensional random equation with cubic term}\label{section:cubicequation}
We consider the equation with a cubic term
\begin{align}\label{cubicequation}
    \frac{\partial u(x,t,\omega)}{\partial t} = D\Delta u(x,t,\omega) - K(\omega)u(x,t,\omega)^3,\quad u(x,0,\omega) = \cos(\pi x),\quad x\in (-1,1),\ t\in (0,T].
\end{align}
The plots in Figure \ref{cubic_D=0} show that in the case $D=0$, the iPCE schemes can produce somewhat competitive errors compared to the niPCE schemes for the EE and ETD-RDP cases, but in the ETDRK4 scheme error stays small more consistently in the niPCE case.\\
For the $D=1$ case, it is seen in Figure \ref{cubic_D=1} that both ETD-RDP and ETDRK4 schemes can produce competitive errors in the iPCE case. However, this has to be put into context taking into account the low stability of iPCE in this case: The step sizes for the iPCE schemes have to be increased substantially in order to produce these errors, resulting in much larger computation times (see Table \ref{parametertable}, case $D=1$, cubic).  
\subsection{Two-dimensional random equation with linear term}
The equation 
\begin{align}\label{2D_linearequation}
    \frac{\partial u(\boldsymbol{x},t,\omega)}{\partial t} = D\Delta u(\boldsymbol{x},t,\omega) - K(\omega)u(\boldsymbol{x},t,\omega),\quad u(\boldsymbol{x},0,\omega) = \cos(\pi x_1)\cos(\pi x_2)
\end{align}
with $\boldsymbol{x} = (x_1,x_2)^\top\in (-1,1)^2,\ t\in (0,T]$, has the exact solution $u(\boldsymbol{x},t,\omega) = \exp(-(K(\omega)+2D\pi^2)t)\cos(\pi x_1)\cos(\pi x_2)$. If $K\sim \mathcal{U}[a,b]$, the expected value is given by
\begin{align}\nonumber
        \mathbb{E}[u(\boldsymbol{x},t,\cdot)] &= \frac{1}{b-a}\int_a^b \exp(-(\xi+2D \pi^2)t)\cos(\pi x_1)\cos(\pi x_2)\; \de \xi \\ \label{2D_exactmean}
        &= \frac{1}{(b-a)t}\cos(\pi x_1)\cos(\pi x_2)  \exp(-(2D\pi^2+ a)t)-\exp(-(2D\pi^2 + b)t).
\end{align}
Figure \ref{D=0_2D} shows the time-dependent errors for the case $D=0$, using the same parameters and time step numbers as in the 1D cases shown in Figure \ref{linear_D=0} and \ref{linear_D=1}. 
\begin{figure}
\newcommand{\X}{2.5}
\newcommand{\Yone}{1.0}
\newcommand{\Ytwo}{0.5}
\newcommand{\intra}{errorarray_FD_EE_mean_system=8_D=0.00000.txt}
\newcommand{\intrb}{errorarray_FD_ETDRDP_mean_system=8_D=0.00000.txt}
\newcommand{\intrc}{errorarray_Spectral_mean_system=8_D=0.00000.txt}
\newcommand{\nintra}{errorarray_mean_Nonintrusive_FD_EE_system=8_D=0.00000.txt}
\newcommand{\nintrb}{errorarray_mean_Nonintrusive_FD_ETDRDP_system=8_D=0.00000.txt}
\newcommand{\nintrc}{errorarray_mean_Nonintrusive_Spectral_system=8_D=0.00000.txt}
\newcommand{\legendposition}{south east}
\centering

\begin{minipage}{\textwidth} 
\centering 
\captionsetup{justification=centering} 
\caption*{\underline{Equation with cubic term, $D=0$, mean}}
\end{minipage}
 \begin{tabular}{@{}c@{}}
\begin{tikzpicture}
    \begin{axis}
    [ cycle list name=goodcolors,
    title={Intrusive PCE, \textbf{EE}},
      height=\Height, width=\Width,
      xlabel={$t$}, ylabel={relative $L^2$ error},ymode = log,
      legend entries={N=1, N=2, N=3, N=4, N=5},
      legend style={at={(0.5,0.475)}, anchor=north} ]
    \foreach \i in {1,2,3,4,5} {
    \addplot table [col sep=comma, x index=0, y index=\i] {\intra}; }
    \end{axis}
\end{tikzpicture}
\begin{tikzpicture}
    \begin{axis}
    [ cycle list name=goodcolors,
    title={Non-Intrusive PCE, \textbf{EE}},
      height=\Height, width=\Width,
      xlabel={$t$}, ylabel={relative $L^2$ error},ymode = log,
      legend entries={MC, QMC, GQ},
      legend pos=south east, ]
    \foreach \i in {1,2,3} {
    \addplot table [col sep=comma, x index=0, y index=\i] {\nintra}; }
    \end{axis}
\end{tikzpicture}
\\
\begin{tikzpicture}
    \begin{axis}
    [ cycle list name=goodcolors,
    title={Intrusive PCE, \textbf{ETD-RDP}},
      height=\Height, width=\Width,
      xlabel={$t$}, ylabel={relative $L^2$ error},ymode = log,
      legend pos=\legendposition, ]
    \foreach \i in {1,2,3,4,5} {
    \addplot table [col sep=comma, x index=0, y index=\i] {\intrb}; }
    \end{axis}
\end{tikzpicture}
\begin{tikzpicture}
    \begin{axis}
    [ cycle list name=goodcolors,
    title={Non-Intrusive PCE, \textbf{ETD-RDP}},
      height=\Height, width=\Width,
      xlabel={$t$}, ylabel={relative $L^2$ error},ymode = log,
      legend pos=\legendposition, ]
    \foreach \i in {1,2,3} {
    \addplot table [col sep=comma, x index=0, y index=\i] {\nintrb}; }
    \end{axis}
\end{tikzpicture}
\\
\begin{tikzpicture}
    \begin{axis}
    [ cycle list name=goodcolors,
    title={Intrusive PCE, \textbf{ETDRK4}},
      height=\Height, width=\Width,
      xlabel={$t$}, ylabel={relative $L^2$ error},ymode = log,
      legend pos=\legendposition, ]
    \foreach \i in {1,2,3,4,5} {
    \addplot table [col sep=comma, x index=0, y index=\i] {\intrc}; }
    \end{axis}
\end{tikzpicture}
\begin{tikzpicture}
    \begin{axis}
    [ cycle list name=goodcolors,
    title={Non-Intrusive PCE, \textbf{ETDRK4}},
      height=\Height, width=\Width,
      xlabel={$t$}, ylabel={relative $L^2$ error},ymode = log,
      legend pos=\legendposition, ]
    \foreach \i in {1,2,3} {
    \addplot table [col sep=comma, x index=0, y index=\i] {\nintrc}; }
    \end{axis}
\end{tikzpicture}
\end{tabular}
\caption{Plots showing the relative time-dependent error for iPCE schemes (Algorithms \ref{algorithm:iPCE_EE}, \ref{algorithm:iPCE_ETDRDP}, \ref{algorithm:iPCE_ETDRK4}, left-hand side) and niPCE schemes (right-hand side) for equation \eqref{cubicequation} with cubic term with diffusion constant $D=0$.}\label{cubic_D=0}
\end{figure}
\begin{figure}
\newcommand{\X}{2.5} \newcommand{\Yone}{1.0}
\newcommand{\Ytwo}{0.5} \newcommand{\intra}{errorarray_FD_EE_mean_system=8_D=1.00000.txt}
\newcommand{\intrb}{errorarray_FD_ETDRDP_mean_system=8_D=1.00000.txt}
\newcommand{\intrc}{errorarray_Spectral_mean_system=8_D=1.00000.txt}
\newcommand{\nintra}{errorarray_mean_Nonintrusive_FD_EE_system=8_D=1.00000.txt}
\newcommand{\nintrb}{errorarray_mean_Nonintrusive_FD_ETDRDP_system=8_D=1.00000.txt}
\newcommand{\nintrc}{errorarray_mean_Nonintrusive_Spectral_system=8_D=1.00000.txt}
\newcommand{\legendposition}{north west}
\centering
\begin{minipage}{\textwidth} 
\centering 
\captionsetup{justification=centering} 
\caption*{\underline{Equation with cubic term, $D=1$, mean}}
\end{minipage}
 \begin{tabular}{@{}c@{}}
\begin{tikzpicture}
    \begin{axis}
    [ cycle list name=goodcolors,
    title={Intrusive PCE, \textbf{EE}},
      height=\Height, width=\Width,
      xlabel={$t$}, ylabel={relative $L^2$ error},ymode = log,
      legend entries={N=1, N=2, N=3, N=4, N=5},
      legend pos=south east, ]
    \foreach \i in {1,2,3,4,5} {
    \addplot table [col sep=comma, x index=0, y index=\i] {\intra}; }
    \end{axis}
\end{tikzpicture}
\begin{tikzpicture}
    \begin{axis}
    [ cycle list name=goodcolors,
    title={Non-Intrusive PCE, \textbf{EE}},
      height=\Height, width=\Width,
      xlabel={$t$}, ylabel={relative $L^2$ error},ymode = log,
      legend entries={MC, QMC, GQ},
      legend pos=south east, ]
    \foreach \i in {1,2,3} {
    \addplot table [col sep=comma, x index=0, y index=\i] {\nintra}; }
    \end{axis}
\end{tikzpicture}
\\
\begin{tikzpicture}
    \begin{axis}
    [ cycle list name=goodcolors,
    title={Intrusive PCE, \textbf{ETD-RDP}},
      height=\Height, width=\Width,
      xlabel={$t$}, ylabel={relative $L^2$ error},ymode = log,
      legend pos=\legendposition, ]
    \foreach \i in {1,2,3,4,5} {
    \addplot table [col sep=comma, x index=0, y index=\i] {\intrb}; }
    \end{axis}
\end{tikzpicture}
\begin{tikzpicture}
    \begin{axis}
    [ cycle list name=goodcolors,
    title={Non-Intrusive PCE, \textbf{ETD-RDP}},
      height=\Height, width=\Width,
      xlabel={$t$}, ylabel={relative $L^2$ error},ymode = log,
      legend pos=\legendposition, ]
    \foreach \i in {1,2,3} {
    \addplot table [col sep=comma, x index=0, y index=\i] {\nintrb}; }
    \end{axis}
\end{tikzpicture}
\\
\begin{tikzpicture}
    \begin{axis}
    [ cycle list name=goodcolors,
    title={Intrusive PCE, \textbf{ETDRK4}},
      height=\Height, width=\Width,
      xlabel={$t$}, ylabel={relative $L^2$ error},ymode = log,
      legend pos=\legendposition, ]
    \foreach \i in {1,2,3,4,5} {
    \addplot table [col sep=comma, x index=0, y index=\i] {\intrc}; }
    \end{axis}
\end{tikzpicture}
\begin{tikzpicture}
    \begin{axis}
    [ cycle list name=goodcolors,
    title={Non-Intrusive PCE, \textbf{ETDRK4}},
      height=\Height, width=\Width,
      xlabel={$t$}, ylabel={relative $L^2$ error},ymode = log,
      legend pos=\legendposition, ]
    \foreach \i in {1,2,3} {
    \addplot table [col sep=comma, x index=0, y index=\i] {\nintrc}; }
    \end{axis}
\end{tikzpicture}
\end{tabular}
\caption{Plots showing the relative time-dependent error for iPCE schemes (Algorithms \ref{algorithm:iPCE_EE}, \ref{algorithm:iPCE_ETDRDP}, \ref{algorithm:iPCE_ETDRK4}, left-hand side) and niPCE schemes (right-hand side) for equation \eqref{cubicequation} with cubic term with diffusion constant $D=1$.}\label{cubic_D=1}
\end{figure}
\begin{figure}
\newcommand{\X}{2.5} \newcommand{\Yone}{1.0}
\newcommand{\Ytwo}{0.5} \newcommand{\faa}{Performanceplot_mean_system=6_D=0.00000.txt}
\newcommand{\fab}{Performanceplot_mean_system=6_D=1.00000.txt}
\newcommand{\fba}{Performanceplot_mean_system=7_D=0.00000.txt}
\newcommand{\fbb}{Performanceplot_mean_system=7_D=1.00000.txt}
\newcommand{\fca}{Performanceplot_mean_system=8_D=0.00000.txt}
\newcommand{\fcb}{Performanceplot_mean_system=8_D=1.00000.txt}
\newcommand{\legendposition}{north west}
\centering
\begin{minipage}{\textwidth} 
\centering 
\captionsetup{justification=centering} 
\caption*{\underline{Performance plots for equations with linear, quadratic and cubic terms}}
\end{minipage}
\begin{tabular}{@{}c@{}}
\begin{tikzpicture}
    \begin{axis}
    [ cycle list name=goodcolors,
    title={Linear, $D=0$},
      height=\Height, width=\Width,
      xlabel={$M$}, ylabel={relative $L^2$ error},ymode = log,xmode = log,
      legend entries={EEi, RDPi, RK4i},
      legend pos=south west, ]
    \foreach \i in {1,2,3,4,5,6} {
    \addplot table [col sep=comma, x index=0, y index=\i] {\faa}; }
    \end{axis}
\end{tikzpicture}
\begin{tikzpicture}
    \begin{axis}
    [ cycle list name=goodcolorsalt,
    title={Linear, $D=1$},
      height=\Height, width=\Width,
      xlabel={$M$}, ylabel={relative $L^2$ error},ymode = log,xmode = log,
      legend entries={EEn, RDPn, RK4n},
      legend style={at={(0.77,0.53)}, anchor=north} ]
      \addlegendimage{color=plotcolor4,densely dashed, dash pattern=on 4pt off 2pt, line width=2pt}
    \addlegendimage{color=plotcolor5,dashed, dash pattern=on 6pt off 2pt, line width=2pt}
    \addlegendimage{color=plotcolor6,style=solid, line width=1pt}
    \foreach \i in {1,2,3,4} {
    \addplot table [col sep=comma, x index=0, y index=\i] {\fab}; }
    \end{axis}
\end{tikzpicture}
\\
\begin{tikzpicture}
    \begin{axis}
    [ cycle list name=goodcolors,
    title={Quadratic, $D=0$},
      height=\Height, width=\Width,
      xlabel={$M$}, ylabel={relative $L^2$ error},ymode = log,xmode = log,
      legend pos=\legendposition, ]
    \foreach \i in {1,2,3,4,5,6} {
    \addplot table [col sep=comma, x index=0, y index=\i] {\fba}; }
    \end{axis}
\end{tikzpicture}
\begin{tikzpicture}
    \begin{axis}
    [ cycle list name=goodcolorsalt,
    title={Quadratic, $D=1$},
      height=\Height, width=\Width,
      xlabel={$M$}, ylabel={relative $L^2$ error},ymode = log,xmode = log,
      legend pos=\legendposition, ]
    \foreach \i in {1,2,3,4} {
    \addplot table [col sep=comma, x index=0, y index=\i] {\fbb}; }
    \end{axis}
\end{tikzpicture}
\\
\begin{tikzpicture}
    \begin{axis}
    [ cycle list name=goodcolors,
    title={Cubic, $D=0$},
      height=\Height, width=\Width,
      xlabel={$M$}, ylabel={relative $L^2$ error},ymode = log,xmode = log,
      legend pos=\legendposition, ]
    \foreach \i in {1,2,3,4,5,6} {
    \addplot table [col sep=comma, x index=0, y index=\i] {\fca}; }
    \end{axis}
\end{tikzpicture}
\begin{tikzpicture}
    \begin{axis}
    [ cycle list name=goodcolorsalt,
    title={Cubic, $D=1$},
      height=\Height, width=\Width,
      xlabel={$M$}, ylabel={relative $L^2$ error},ymode = log,xmode = log,
      legend pos=\legendposition, ]
    \foreach \i in {1,2,3,4} {
    \addplot table [col sep=comma, x index=0, y index=\i] {\fcb}; }
    \end{axis}
\end{tikzpicture}
\end{tabular}
\caption{Plots showing the error for the different equations at time $T=2$ (except for equation \ref{quadraticequation} with the quadratic term and $D=0$, where $T=0.4$). The shorthands EEi, EEn in the legend stand for EE iPCE and EE niPCE, respectively, and analogously with RDP for ETD-RDP and RK4 for ETDRK4. For $D=1$, only ETD-RDP and ETDRK4 are shown, since EE requires higher $M$ in order to be stable.}\label{performanceplots}
\end{figure}
\subsection{Performance plots}\label{section:performance_plots}
In Figure \ref{performanceplots}, for the different equations we show how the error for iPCE and niPCE compare for different numbers of time steps $M$. For the niPCE errors, $q=10$ Gauss-Legendre quadrature points were used. It can be seen that in most cases, the iPCE and niPCE errors are very similar up to a certain $M$, where the error lines split up and iPCE does not go below a certain threshold. This threshold can, in some cases, be higher than the niPCE error, so that the iPCE error appears as a constant line, independent of $M$ (as in the quadratic ad cubic case for $D=0$). This additional error for iPCE may be caused by the approximation of the iPCE product terms given in \eqref{u3expansion}.
\subsection{Random Gray-Scott model}\label{section:grayscott}
We consider the random Gray-Scott model
\begin{align}\label{grayscottsystem}
    \begin{cases}
        \frac{\partial u(x,t,\omega)}{\partial t} = D_u \Delta u(x,t,\omega) -u(x,t,\omega)v(x,t,\omega)^2 + F(1-u(x,t,\omega)), \\
        \frac{\partial v(x,t,\omega)}{\partial t} = D_v \Delta v(x,t,\omega) + u(x,t,\omega)v(x,t,\omega)^2 - (F+k(\omega))v(x,t,\omega).
    \end{cases}
\end{align}
for $x\in (-1,1)$, $t\in (0,T]$ for some $T>0$, $D_u=2\cdot 10^{-5}$, $D_v=10^{-5}$, a constant $F$ and a uniformly distributed random variable $k$. We take the initial condition from \cite[Figure 9]{xxpaut} which is given by
\begin{align*}
    u_{\mathrm{init}}(x) = 1-\frac{5}{3\sqrt{2\pi}}\exp\left( -6(x-\mu)^2\right),\quad 
    v_{\mathrm{init}}(x) = 0.37\cdot\frac{7.5}{2\sqrt{2}\Gamma(\frac{1}{3})}\exp\left( \frac{-7(x-\mu)^3}{\sqrt{2}} \right)
\end{align*}
where $\mu = 0$ is the midpoint of the interval. For the two-dimensional plots shown in Figures \ref{figure:grayscottsensitivity} and \ref{brokenpattern}, we use as an initial condition a function with four off-center local extrema:
\begin{align}\label{2dinitialcondition}
    v(x,y,0) = \frac{1}{4}\sum_{i=1}^4 \exp\left(-150((x-x_i)^2+(y-y_i)^2)\right),\quad u(x,y,0) = 1-v(x,y,0),
\end{align}
where $(x_1,y_1) = (2/7,2/7)$, $(x_2,y_2) = (-2/7,2/7)$, $(x_3,y_3) = (2/7,-2/7)$, $(x_4,y_4) = (-2/7,-2/7)$.
\begin{remark}
    The Gray-Scott model has been extensively studied and is known to show a wide variety of pattern formation behavior. The patterns are called Turing patterns, for an overview we refer to \cite{pearson1993complex,grayscottonline}. Depending on the parameter values of $k$ and $F$, $u$ and $v$ may enter trivial or non-trivial homogeneous steady states, or non-homogeneous states called Turing patterns. A necessary condition in the $k$-$F$ parameter space in order for Turing pattern formation to occur is \cite[Eq. (19)]{mazin1996pattern}
\begin{align}\label{turingregion}
    [2(F+k)-(v_0^2+F)]^2>8(F+k)(v_0^2-F),
\end{align} 
where $v_0$ is either the trivial steady state (the `red state') $v_R=0$ or, if $d:= 1-4(F+k)^2/F>0$ holds, one of the two non-trivial steady states \cite[Eqs. (7), (8)]{mazin1996pattern} (with $v_B$ also called the `blue state')
\begin{align*}
    v_B = \frac{1}{2}\alpha(1+\sqrt{d}),\quad v_1 = \frac{1}{2}\alpha(1-\sqrt{d}).
\end{align*}
\end{remark}
Since now there are two coupled functions instead of one, the schemes need to take into account these two functions, but for the sake of brevity we will not spell out all the schemes again, but refer to \cite{Kleefeld2020,kassamsolving} (the intrusive PCE schemes given in Section \ref{section:intrusivepce} can then easily be extended to two functions).\\
We pick for all simulations (1D and 2D) $F=0.04$ and $k(\omega) \sim \mathcal{U}[0.058,0.062]$, which falls into the region of complex pattern formation described by \eqref{turingregion}. It is seen in Figure \ref{grayscotterror} that the iPCE schemes fail to compute correct solutions for the Gray-Scott model. Furthermore, the needed time step numbers in order for the iPCE schemes are even higher than for Equation \ref{cubicequation} with a single cubic term. This leaves niPCE as the only viable option to treat the random Gray-Scott system. We show an example of such a simulation in two spatial dimensions for a mean $\mathbb{E}[u(x,t,\cdot)]$ in Figure \ref{brokenpattern} with the correct result produced by niPCE juxtaposed by iPCE simulations which fail to reproduce the pattern. This is likely due to the challenge of sharp dependency in the random variable, as in the seemingly small interval $k\in[0.058,0.062]$, a wide range of patterns emerge (see also \cite{grayscottonline}). Furthermore, long-term integration of the Gray-Scott system is a challenge even in the deterministic case, as show in Figure \ref{figure:grayscottsensitivity} for $T=5000$, where it is seen that a high level of accuracy in both space and time is needed in order to produce the correct pattern shown on the bottom right. Especially in the cases of EE and ETD-RDP-IF, for $p=256$ the error of the second-order spatial finite difference approximation appears to dominate, causing the pattern to look rather different.\\
Given the difficulties of long-term integration for the Gray-Scott model and the numerical challenges and substantially increased computation time for iPCE, it appears to be advisable to use niPCE for systems with strong nonlinearities such as Gray-Scott.
\begin{figure}[h]
\newcommand{\X}{2.5}
\newcommand{\Yone}{1.0}
\newcommand{\Ytwo}{0.5}
\newcommand{\intra}{errorarray_FD_EE_2D_mean_system=6_D=1.00000.txt}
\newcommand{\intrb}{errorarray_FD_ETDRDP_2D_mean_system=6_D=1.00000.txt}
\newcommand{\intrc}{errorarray_Spectral_2D_mean_system=6_D=1.00000.txt}
\newcommand{\legendposition}{north west}
\centering
\begin{minipage}{\textwidth}
\centering 
\captionsetup{justification=centering} 
\caption*{\underline{iPCE errors for 2D random equation with linear term, $D=1$}}
\end{minipage}
 \begin{tabular}{@{}c@{}}
\begin{tikzpicture}
    \begin{axis}
    [ cycle list name=goodcolors,
    title={\textbf{EE}},
      height=\Heightt, width=\Widtht,
      xlabel={$t$}, ylabel={relative $L^2$ error},ymode = log,
      legend entries={$N=1$,$N=2$,$N=3$,$N=4$,$N=5$},
      legend pos=south east, ]
    \foreach \i in {1,2,3,4,5} {
    \addplot table [col sep=comma, x index=0, y index=\i] {\intra};
    }
    \end{axis}
\end{tikzpicture}
\begin{tikzpicture}
    \begin{axis}
    [ cycle list name=goodcolors,
    title={\textbf{ETD-RDP-IF}},
      height=\Heightt, width=\Widtht,
      xlabel={$t$},ymode = log,
      legend pos=north west, ]
 \addlegendimage{plotcolor3,dotted, line width=2pt}
        \addlegendimage{plotcolor4,densely dashed, line width=2pt}
    \foreach \i in {1,2,3,4,5} {
    \addplot table [col sep=comma, x index=0, y index=\i] {\intrb};
    }
    \end{axis}
\end{tikzpicture}
\begin{tikzpicture}
    \begin{axis}
    [ cycle list name=goodcolors,
    title={\textbf{ETDRK4}},
      height=\Heightt, width=\Widtht,
      xlabel={$t$},ymode = log,
      legend pos=north west, ]
      \addlegendimage{plotcolor5,dashed, dash pattern=on 6pt off 2pt, line width=2pt}
    \foreach \i in {1,2,3,4,5} {
    \addplot table [col sep=comma, x index=0, y index=\i] {\intrc};
    }
    \end{axis}
\end{tikzpicture}
\end{tabular}
\caption{Time-dependent $L^2$ errors for equation \eqref{2D_linearequation}, using the iPCE algorithms \ref{algorithm:iPCE_EE}, \ref{algorithm:iPCE_ETDRDPIF} and \ref{algorithm:iPCE_ETDRK4}. We used the same parameters and numbers of time steps as in the 1D case from Figure \ref{linear_D=0} and Table \ref{parametertable}. The exact expected value is given by \eqref{2D_exactmean}. The spatial resolution in the ETD-RDP case was picked in this simulation as $p=512$. For lower $p$, the error caused by the finite difference discretization dominates and the curves for different $N$ are identical. } \label{D=0_2D}
\end{figure}
\begin{figure}[h]
\newcommand{\X}{2.5}
\newcommand{\Yone}{1.0}
\newcommand{\Ytwo}{0.5}
\newcommand{\intra}{errorarray_FD_EE_mean_system=0_D=0.00002.txt}
\newcommand{\intrb}{errorarray_FD_ETDRDP_mean_system=0_D=0.00002.txt}
\newcommand{\intrc}{errorarray_Spectral_mean_system=0_D=0.00002.txt}
\newcommand{\legendposition}{north west}
\centering
\begin{minipage}{\textwidth} 
\centering 
\captionsetup{justification=centering} 
\caption*{\underline{Errors for a random Gray-Scott system}}
\end{minipage}
 \begin{tabular}{@{}c@{}}
\begin{tikzpicture}
    \begin{axis}
    [ cycle list name=goodcolors,
    title={\textbf{EE}},
      height=\Heightt, width=\Widtht,
      xlabel={$t$}, ylabel={relative $L^2$ error},ymode = log,
      legend entries={$N=1$,$N=2$,$N=3$,$N=4$,$N=5$},
      legend pos=south east, ]
    \foreach \i in {1,2,3,4,5} {
    \addplot table [col sep=comma, x index=0, y index=\i] {\intra};
    }
    \end{axis}
\end{tikzpicture}
\begin{tikzpicture}
    \begin{axis}
    [ cycle list name=goodcolors,
    title={\textbf{ETD-RDP}},
      height=\Heightt, width=\Widtht,
      xlabel={$t$},ymode = log,
      legend pos=north west, ]
 \addlegendimage{plotcolor3,dotted, line width=2pt}
        \addlegendimage{plotcolor4,densely dashed, line width=2pt}
    \foreach \i in {1,2,3,4,5} {
    \addplot table [col sep=comma, x index=0, y index=\i] {\intrb};
    }
    \end{axis}
\end{tikzpicture}
\begin{tikzpicture}
    \begin{axis}
    [ cycle list name=goodcolors,
    title={\textbf{ETDRK4}},
      height=\Heightt, width=\Widtht,
      xlabel={$t$},ymode = log,
      legend pos=north west, ]
      \addlegendimage{plotcolor5,dashed, dash pattern=on 6pt off 2pt, line width=2pt}
    \foreach \i in {1,2,3,4,5} {
    \addplot table [col sep=comma, x index=0, y index=\i] {\intrc};
    }
    \end{axis}
\end{tikzpicture}
\end{tabular}
\caption{Time-dependent error plot for a random 1D Gray-Scott system. None of the iPCE schemes work in this case: Initially, the solution is relatively accurate, but as pattern formation starts to occur, iPCE breaks down. For a visual example in two dimensions, see Figure \ref{brokenpattern}.} \label{grayscotterror}
\end{figure}
\begin{figure}[h]
\newcommand{\aax}{0}
\newcommand{\aay}{0}
\newcommand{\bbx}{256}
\newcommand{\bby}{0}
\newcommand{\ccx}{512}
\newcommand{\ccy}{0}
\newcommand{\ddx}{0}
\newcommand{\ddy}{-270}
\newcommand{\eex}{256}
\newcommand{\eey}{-270}
\newcommand{\ffx}{512}
\newcommand{\ffy}{-270}
\newcommand{\mycolormap}{colormap={mymap}{[1pt] rgb(0pt)=(0.2422,0.1504,0.6603); rgb(1pt)=(0.25039,0.164995,0.707614); rgb(2pt)=(0.257771,0.181781,0.751138); rgb(3pt)=(0.264729,0.197757,0.795214); rgb(4pt)=(0.270648,0.214676,0.836371); rgb(5pt)=(0.275114,0.234238,0.870986); rgb(6pt)=(0.2783,0.255871,0.899071); rgb(7pt)=(0.280333,0.278233,0.9221); rgb(8pt)=(0.281338,0.300595,0.941376); rgb(9pt)=(0.281014,0.322757,0.957886); rgb(10pt)=(0.279467,0.344671,0.971676); rgb(11pt)=(0.275971,0.366681,0.982905); rgb(12pt)=(0.269914,0.3892,0.9906); rgb(13pt)=(0.260243,0.412329,0.995157); rgb(14pt)=(0.244033,0.435833,0.998833); rgb(15pt)=(0.220643,0.460257,0.997286); rgb(16pt)=(0.196333,0.484719,0.989152); rgb(17pt)=(0.183405,0.507371,0.979795); rgb(18pt)=(0.178643,0.528857,0.968157); rgb(19pt)=(0.176438,0.549905,0.952019); rgb(20pt)=(0.168743,0.570262,0.935871); rgb(21pt)=(0.154,0.5902,0.9218); rgb(22pt)=(0.146029,0.609119,0.907857); rgb(23pt)=(0.138024,0.627629,0.89729); rgb(24pt)=(0.124814,0.645929,0.888343); rgb(25pt)=(0.111252,0.6635,0.876314); rgb(26pt)=(0.0952095,0.679829,0.859781); rgb(27pt)=(0.0688714,0.694771,0.839357); rgb(28pt)=(0.0296667,0.708167,0.816333); rgb(29pt)=(0.00357143,0.720267,0.7917); rgb(30pt)=(0.00665714,0.731214,0.766014); rgb(31pt)=(0.0433286,0.741095,0.73941); rgb(32pt)=(0.0963952,0.75,0.712038); rgb(33pt)=(0.140771,0.7584,0.684157); rgb(34pt)=(0.1717,0.766962,0.655443); rgb(35pt)=(0.193767,0.775767,0.6251); rgb(36pt)=(0.216086,0.7843,0.5923); rgb(37pt)=(0.246957,0.791795,0.556743); rgb(38pt)=(0.290614,0.79729,0.518829); rgb(39pt)=(0.340643,0.8008,0.478857); rgb(40pt)=(0.3909,0.802871,0.435448); rgb(41pt)=(0.445629,0.802419,0.390919); rgb(42pt)=(0.5044,0.7993,0.348); rgb(43pt)=(0.561562,0.794233,0.304481); rgb(44pt)=(0.617395,0.787619,0.261238); rgb(45pt)=(0.671986,0.779271,0.2227); rgb(46pt)=(0.7242,0.769843,0.191029); rgb(47pt)=(0.773833,0.759805,0.16461); rgb(48pt)=(0.820314,0.749814,0.153529); rgb(49pt)=(0.863433,0.7406,0.159633); rgb(50pt)=(0.903543,0.733029,0.177414); rgb(51pt)=(0.939257,0.728786,0.209957); rgb(52pt)=(0.972757,0.729771,0.239443); rgb(53pt)=(0.995648,0.743371,0.237148); rgb(54pt)=(0.996986,0.765857,0.219943); rgb(55pt)=(0.995205,0.789252,0.202762); rgb(56pt)=(0.9892,0.813567,0.188533); rgb(57pt)=(0.978629,0.838629,0.176557); rgb(58pt)=(0.967648,0.8639,0.16429); rgb(59pt)=(0.96101,0.889019,0.153676); rgb(60pt)=(0.959671,0.913457,0.142257); rgb(61pt)=(0.962795,0.937338,0.12651); rgb(62pt)=(0.969114,0.960629,0.106362); rgb(63pt)=(0.9769,0.9839,0.0805)},
}
\begin{tikzpicture}
\begin{axis}[%
width=\Widthgs,
height=\Heightgs,
at={(\aax,\aay)},
scale only axis,
point meta min=0.000215671688202298,
point meta max=0.378861626873836,
axis on top,
xmin=-1,
xmax=1,
ymin=-1,
ymax=1,
axis background/.style={fill=white},
title style={font=\bfseries},
title={EE, $p=256$},
\mycolormap,
]
\addplot [forget plot] graphics [xmin=-1.00390625, xmax=0.99609375, ymin=-1.00390625, ymax=0.99609375] {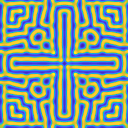};
\end{axis}
\begin{axis}[%
width=\Widthgs,
height=\Heightgs,
at={(\bbx,\bby)},
scale only axis,
point meta min=0.00100986762483131,
point meta max=0.376638920408934,
axis on top,
xmin=-1,
xmax=1,
ymin=-1,
ymax=1,
axis background/.style={fill=white},
title style={font=\bfseries},
title={ETD-RDP-IF, $p=256$},
\mycolormap,
]
\addplot [forget plot] graphics [xmin=-1.0009765625, xmax=0.9990234375, ymin=-1.0009765625, ymax=0.9990234375] {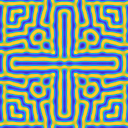};
\end{axis}
\begin{axis}[%
width=\Widthgs,
height=\Heightgs,
at={(\ccx,\ccy)},
scale only axis,
point meta min=0.00100986762483131,
point meta max=0.376638920408934,
axis on top,
xmin=-1,
xmax=1,
ymin=-1,
ymax=1,
axis background/.style={fill=white},
title style={font=\bfseries},
title={ETDRK4, $p=256$},
\mycolormap,
colorbar
]
\addplot [forget plot] graphics [xmin=-1.0009765625, xmax=0.9990234375, ymin=-1.0009765625, ymax=0.9990234375] {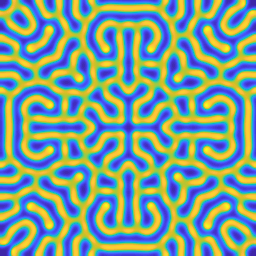};
\end{axis}\bigskip
\begin{axis}[%
width=\Widthgs,
height=\Heightgs,
at={(\ddx,\ddy)},
scale only axis,
point meta min=0.00100986762483131,
point meta max=0.376638920408934,
axis on top,
xmin=-1,
xmax=1,
ymin=-1,
ymax=1,
axis background/.style={fill=white},
title style={font=\bfseries},
title={EE, $p=1024$},
\mycolormap,
]
\addplot [forget plot] graphics [xmin=-1.0009765625, xmax=0.9990234375, ymin=-1.0009765625, ymax=0.9990234375] {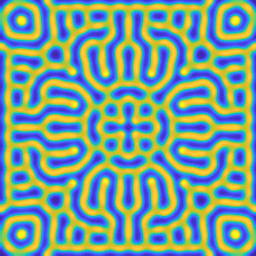};
\end{axis}
\begin{axis}[%
width=\Widthgs,
height=\Heightgs,
at={(\eex,\eey)},
scale only axis,
point meta min=0.00100986762483131,
point meta max=0.376638920408934,
axis on top,
xmin=-1,
xmax=1,
ymin=-1,
ymax=1,
axis background/.style={fill=white},
title style={font=\bfseries},
title={ETD-RDP-IF, $p=1024$},
\mycolormap,
]
\addplot [forget plot] graphics [xmin=-1.0009765625, xmax=0.9990234375, ymin=-1.0009765625, ymax=0.9990234375] {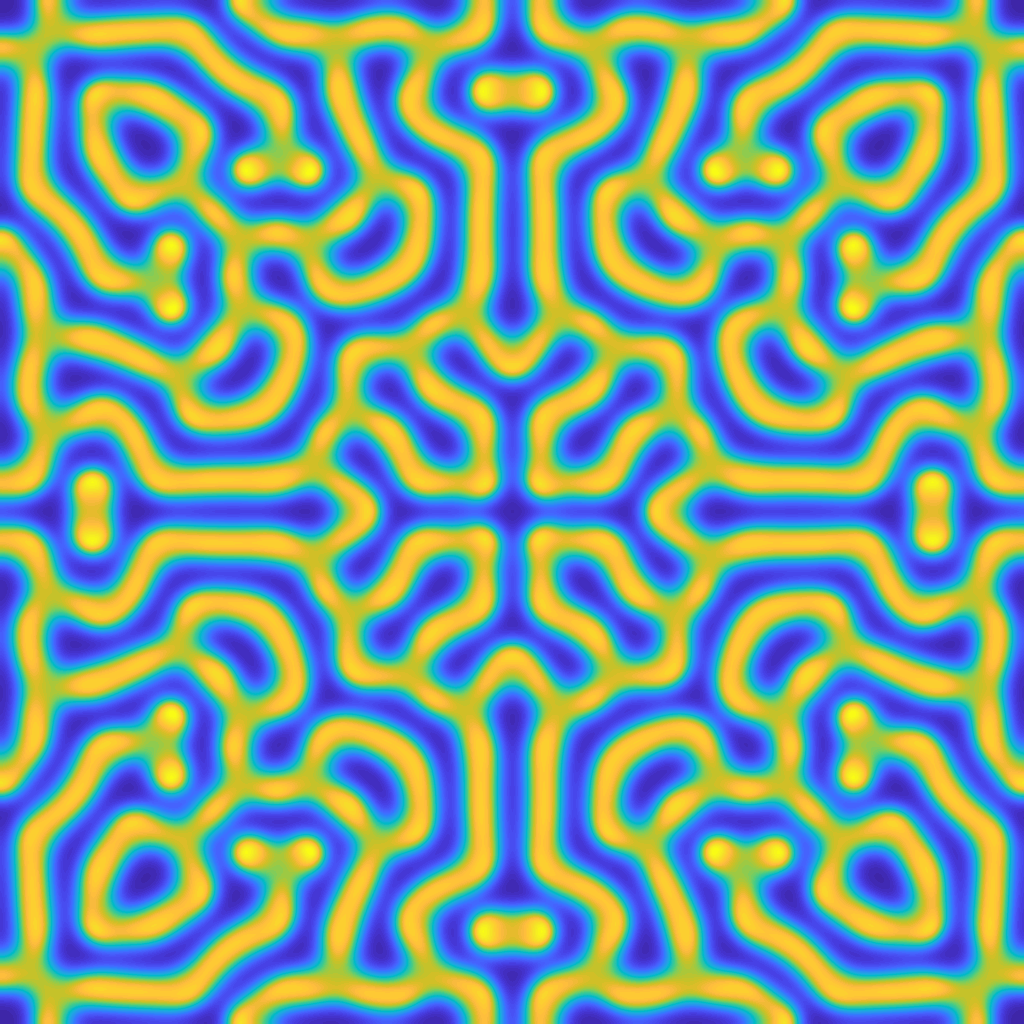};
\end{axis}
\begin{axis}[%
width=\Widthgs,
height=\Heightgs,
at={(\ffx,\ffy)},
scale only axis,
point meta min=0.00100986762483131,
point meta max=0.376638920408934,
axis on top,
xmin=-1,
xmax=1,
ymin=-1,
ymax=1,
axis background/.style={fill=white},
title style={font=\bfseries},
title={ETDRK4, $p=1024$},
\mycolormap,
colorbar
]
\addplot [forget plot] graphics [xmin=-1.0009765625, xmax=0.9990234375, ymin=-1.0009765625, ymax=0.9990234375] {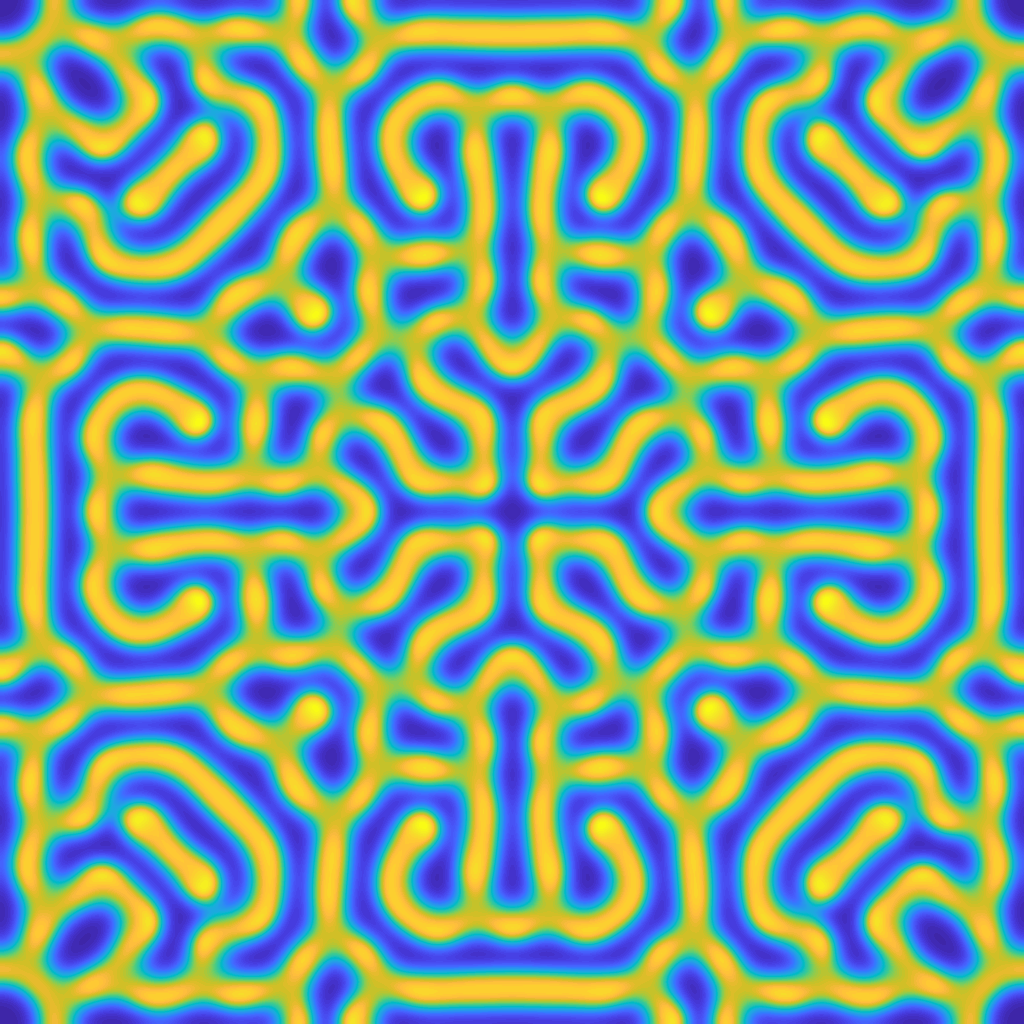};
\end{axis}
\end{tikzpicture}%
\caption{Deterministic simulations of Gray-Scott patterns for different schemes and resolutions for the same initial condition at $T=5000$ (using Algorithms \ref{algorithm:EE_deterministic}, \ref{algorithm:ETDRDPIF_deterministic}, and \ref{algorithm:ETDRK4_deterministic}) . It becomes apparent that after a long-term integration, Gray-Scott patterns are extremely sensitive to the spatial resolution (and number of time steps), and for the correct solution, a high spatial resolution is needed, along with a scheme such as ETDRK4 which performs a spectral approximation in space. The numbers of time steps used are $M=100000$ for EE and $M=10000$ for ETD-RDP-IF and ETDRK4.}\label{figure:grayscottsensitivity}
\end{figure}
\begin{figure}[h]
\newcommand{\aax}{0}
\newcommand{\aay}{0}
\newcommand{\bbx}{256}
\newcommand{\bby}{0}
\newcommand{\ccx}{512}
\newcommand{\ccy}{0}
\newcommand{\ddx}{0}
\newcommand{\ddy}{-270}
\newcommand{\eex}{256}
\newcommand{\eey}{-270}
\newcommand{\ffx}{512}
\newcommand{\ffy}{-270}
\newcommand{\mycolormap}{colormap={mymap}{[1pt] rgb(0pt)=(0.2422,0.1504,0.6603); rgb(1pt)=(0.25039,0.164995,0.707614); rgb(2pt)=(0.257771,0.181781,0.751138); rgb(3pt)=(0.264729,0.197757,0.795214); rgb(4pt)=(0.270648,0.214676,0.836371); rgb(5pt)=(0.275114,0.234238,0.870986); rgb(6pt)=(0.2783,0.255871,0.899071); rgb(7pt)=(0.280333,0.278233,0.9221); rgb(8pt)=(0.281338,0.300595,0.941376); rgb(9pt)=(0.281014,0.322757,0.957886); rgb(10pt)=(0.279467,0.344671,0.971676); rgb(11pt)=(0.275971,0.366681,0.982905); rgb(12pt)=(0.269914,0.3892,0.9906); rgb(13pt)=(0.260243,0.412329,0.995157); rgb(14pt)=(0.244033,0.435833,0.998833); rgb(15pt)=(0.220643,0.460257,0.997286); rgb(16pt)=(0.196333,0.484719,0.989152); rgb(17pt)=(0.183405,0.507371,0.979795); rgb(18pt)=(0.178643,0.528857,0.968157); rgb(19pt)=(0.176438,0.549905,0.952019); rgb(20pt)=(0.168743,0.570262,0.935871); rgb(21pt)=(0.154,0.5902,0.9218); rgb(22pt)=(0.146029,0.609119,0.907857); rgb(23pt)=(0.138024,0.627629,0.89729); rgb(24pt)=(0.124814,0.645929,0.888343); rgb(25pt)=(0.111252,0.6635,0.876314); rgb(26pt)=(0.0952095,0.679829,0.859781); rgb(27pt)=(0.0688714,0.694771,0.839357); rgb(28pt)=(0.0296667,0.708167,0.816333); rgb(29pt)=(0.00357143,0.720267,0.7917); rgb(30pt)=(0.00665714,0.731214,0.766014); rgb(31pt)=(0.0433286,0.741095,0.73941); rgb(32pt)=(0.0963952,0.75,0.712038); rgb(33pt)=(0.140771,0.7584,0.684157); rgb(34pt)=(0.1717,0.766962,0.655443); rgb(35pt)=(0.193767,0.775767,0.6251); rgb(36pt)=(0.216086,0.7843,0.5923); rgb(37pt)=(0.246957,0.791795,0.556743); rgb(38pt)=(0.290614,0.79729,0.518829); rgb(39pt)=(0.340643,0.8008,0.478857); rgb(40pt)=(0.3909,0.802871,0.435448); rgb(41pt)=(0.445629,0.802419,0.390919); rgb(42pt)=(0.5044,0.7993,0.348); rgb(43pt)=(0.561562,0.794233,0.304481); rgb(44pt)=(0.617395,0.787619,0.261238); rgb(45pt)=(0.671986,0.779271,0.2227); rgb(46pt)=(0.7242,0.769843,0.191029); rgb(47pt)=(0.773833,0.759805,0.16461); rgb(48pt)=(0.820314,0.749814,0.153529); rgb(49pt)=(0.863433,0.7406,0.159633); rgb(50pt)=(0.903543,0.733029,0.177414); rgb(51pt)=(0.939257,0.728786,0.209957); rgb(52pt)=(0.972757,0.729771,0.239443); rgb(53pt)=(0.995648,0.743371,0.237148); rgb(54pt)=(0.996986,0.765857,0.219943); rgb(55pt)=(0.995205,0.789252,0.202762); rgb(56pt)=(0.9892,0.813567,0.188533); rgb(57pt)=(0.978629,0.838629,0.176557); rgb(58pt)=(0.967648,0.8639,0.16429); rgb(59pt)=(0.96101,0.889019,0.153676); rgb(60pt)=(0.959671,0.913457,0.142257); rgb(61pt)=(0.962795,0.937338,0.12651); rgb(62pt)=(0.969114,0.960629,0.106362); rgb(63pt)=(0.9769,0.9839,0.0805)},
}
\begin{tikzpicture}
\begin{axis}[%
width=\Widthgs,
height=\Heightgs,
at={(\aax,\aay)},
scale only axis,
point meta min=0.000215671688202298,
point meta max=0.378861626873836,
axis on top,
xmin=-1,
xmax=1,
ymin=-1,
ymax=1,
axis background/.style={fill=white},
title style={font=\bfseries},
title={niPCE ETDRK4, $T=1500$},
\mycolormap,
]
\addplot [forget plot] graphics [xmin=-1.00390625, xmax=0.99609375, ymin=-1.00390625, ymax=0.99609375] {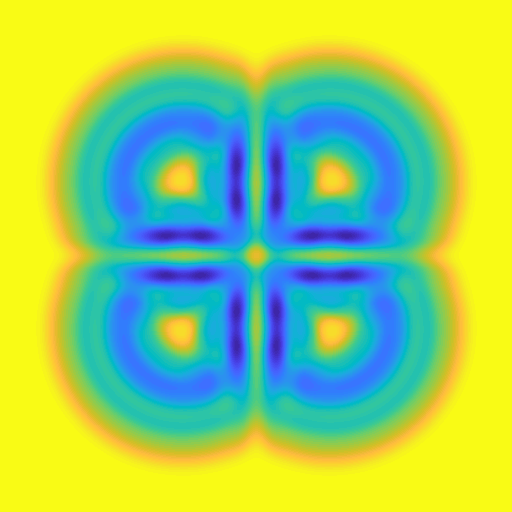};
\end{axis}
\begin{axis}[%
width=\Widthgs,
height=\Heightgs,
at={(\bbx,\bby)},
scale only axis,
point meta min=0.00100986762483131,
point meta max=0.376638920408934,
axis on top,
xmin=-1,
xmax=1,
ymin=-1,
ymax=1,
axis background/.style={fill=white},
title style={font=\bfseries},
title={iPCE EE, $T=500$},
\mycolormap,
]
\addplot [forget plot] graphics [xmin=-1.0009765625, xmax=0.9990234375, ymin=-1.0009765625, ymax=0.9990234375] {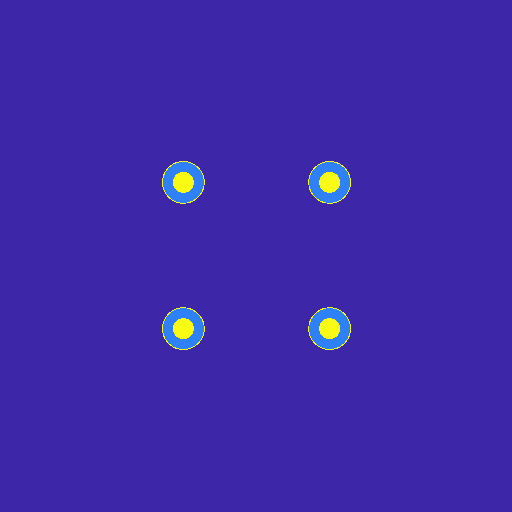};
\end{axis}
\begin{axis}[%
width=\Widthgs,
height=\Heightgs,
at={(\ccx,\ccy)},
scale only axis,
point meta min=0.00100986762483131,
point meta max=0.376638920408934,
axis on top,
xmin=-1,
xmax=1,
ymin=-1,
ymax=1,
axis background/.style={fill=white},
title style={font=\bfseries},
title={iPCE ETDRK4, $T=200$},
\mycolormap,
colorbar
]
\addplot [forget plot] graphics [xmin=-1.0009765625, xmax=0.9990234375, ymin=-1.0009765625, ymax=0.9990234375] {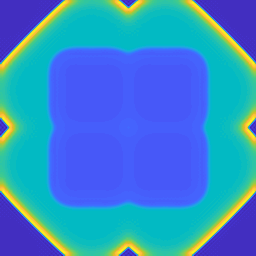};
\end{axis}\bigskip
\end{tikzpicture}%
\caption{Simulations for $\mathbb{E}[u(x,T,\cdot)]$ of the random Gray-Scott system \eqref{grayscottsystem} with $F=0.04$ and $k(\omega)\sim \mathcal{U}[0.058,0.062]$. A niPCE simulation with ETDRK4 (left-hand side) and GQ gives an impression of superimposed patterns for different $k$. Over time, the four off-center bumps of the initial condition \eqref{2dinitialcondition} expand and connect, forming intricate patterns as observed on the left-hand side. The EE iPCE simulation fails to correctly propagate the patterns, and the four rings from the initial condition stop expanding. For the ETDRK4 iPCE simulation, the edges of the rings propagate too fast, and no pattern formation is observed. }\label{brokenpattern}
\end{figure}
\newpage
\section{Conclusion and outlook.}\label{section:Outlook}
In this work, we investigated how the ETD-RDP-IF and ETDRK4 schemes can be implemented in an iPCE scheme and compared the performance of these two schemes and an EE scheme to a niPCE approach. While niPCE using Monte Carlo methods or Gaussian quadrature is in most cases superior to iPCE, we also found that in some cases such as the model equation with a quadratic term, iPCE results in lower errors for comparable runtimes. For complex pattern formation dynamics such as the Gray-Scott model, iPCE breaks down for all three schemes which leaves only the non-intrusive variance as a viable option. While not implemented in full detail, it is also apparent that the curse of dimensionality poses a bigger problem to iPCE than it does to niPCE since solving bigger iPCE systems scales much worse than using sparse grids for niPCE.\\
Future work could include methods to remedy the shortcomings of iPCE for complex dynamical systems. In this work, it was seen that a direct iPCE implementation of algorithms which are powerful in the deterministic case is not sufficient. Works such as \cite{bonnaire2021intrusive} have introduced an asynchronous time integration method for iPCE to deal with sharp dependencies in the random variable. Recent work \cite{eckert2020polynomial} also shows that B-splines can be used instead of classical orthogonal polynomial bases to further reduce the error. One possibility for achieving speedups for the ETD-RDP-IF scheme is the use of parallelization as has been demonstrated in \cite[pp. 8]{Kleefeld2020}. It could be carried out in an analogous fashion for the iPCE scheme.\\
Another topic which could be investigated are non-polynomial nonlinear functions, which can be handled using truncations of Taylor series \cite{debusschere2004}.
\begin{table}[h]
    \centering
    \begin{tabular}{c|c|c}
    \toprule
         \multicolumn{3}{c}{$D=0$, linear} \\
         \midrule
        & Intrusive & Non-Intrusive \\
         \midrule
        \textbf{EE} & 0.3225   & 0.6717  \\
         \textbf{ETD-RDP} & 0.4112  & 0.2256  \\
         \textbf{ETDRK4} & 0.4322 & 0.4244 \\
         \bottomrule
    \end{tabular}
    \qquad\quad
    \begin{tabular}{c|c|c}
    \toprule
         \multicolumn{3}{c}{$D=1$, linear} \\
         \midrule
        & Intrusive & Non-Intrusive \\
         \midrule
        \textbf{EE} & 0.4386  & 6.4825 \\
         \textbf{ETD-RDP} & 0.2741 & 0.3698 \\
         \textbf{ETDRK4} & 0.7304 & 0.4809 \\
         \bottomrule
    \end{tabular}

    \begin{tabular}{c|c|c}

         \multicolumn{3}{c}{$D=0$, quadratic} \\
         \midrule
        & Intrusive & Non-Intrusive \\
         \midrule
        \textbf{EE} &  0.5090  & 0.7563  \\
         \textbf{ETD-RDP} & 0.3143 & 0.2550  \\
         \textbf{ETDRK4} & 0.7242 & 0.6093 \\
         \bottomrule
    \end{tabular}
    \qquad\quad
    \begin{tabular}{c|c|c}
         \multicolumn{3}{c}{$D=1$, quadratic} \\
         \midrule
        & Intrusive & Non-Intrusive \\
         \midrule
        \textbf{EE} &  3.5583  & 5.8424  \\
         \textbf{ETD-RDP} & 0.4906 & 0.3265  \\
         \textbf{ETDRK4} & 1.1047 & 0.4916 \\
         \bottomrule
    \end{tabular}

    \vspace{0.1cm}
    \begin{tabular}{c|c|c}

         \multicolumn{3}{c}{$D=0$, cubic} \\
         \midrule
        & Intrusive & Non-Intrusive \\
         \midrule
        \textbf{EE} &  12.6065  & 0.4934  \\
         \textbf{ETD-RDP} & 5.2208 & 0.4339  \\
         \textbf{ETDRK4} & 6.1375 & 0.5221 \\
         \bottomrule
    \end{tabular}
    \qquad\quad
    \begin{tabular}{c|c|c}

         \multicolumn{3}{c}{$D=1$, cubic} \\
         \midrule
        & Intrusive & Non-Intrusive \\
         \midrule
        \textbf{EE} &  247.6447  & 0.3319  \\
         \textbf{ETD-RDP} & 9.8553 & 0.5722  \\
         \textbf{ETDRK4} & 12.4297 & 0.5626 \\
         \bottomrule
    \end{tabular}
    \caption{Runtimes for the created plots, all times in seconds, for iPCE with $N=5$ and for niPCE with 50 realizations}
    \label{runtimetable}
\end{table}

\begin{table}[h]
\centering
    \begin{tabular}{c|c|cc}
        \multicolumn{4}{c}{$D=0$, linear} \\
         \midrule
         & iPCE & \multicolumn{2}{c}{niPCE} \\
         & $M$ & $M$ & $q$ \\
         \midrule
         \textbf{EE} & 1000 & 2000 & 50\\
         \textbf{ETD-RDP} & 200 & 200 & 50\\
         \textbf{ETDRK4} & 100 & 100 & 50\\
         \textbf{ETDRK4} ref. & & - & -  \\
         \bottomrule
    \end{tabular}
\quad 
    \begin{tabular}{c|c|cc}
        \multicolumn{4}{c}{$D=1$, linear} \\
         \midrule
         & iPCE & \multicolumn{2}{c}{niPCE} \\
         & $M$ & $M$ & $q$ \\
         \midrule
         \textbf{EE} & 20000 & 20000 & 50\\
         \textbf{ETD-RDP} & 400 & 200 & 50\\
         \textbf{ETDRK4} & 200 & 100 & 50\\
         \textbf{ETDRK4} ref. & & 1000 & 200  \\
         \bottomrule
    \end{tabular}\\
    \begin{tabular}{c|c|cc}
        \multicolumn{4}{c}{$D=0$, quadratic} \\
         \midrule
         & iPCE & \multicolumn{2}{c}{niPCE} \\
         & $M$ & $M$ & $q$ \\
         \midrule
         \textbf{EE} & 1000 & 2000 & 50\\
         \textbf{ETD-RDP} & 200 & 200 & 50\\
         \textbf{ETDRK4} & 100 & 100 & 50\\
         \textbf{ETDRK4} ref. & & 1000 & 200  \\
         \bottomrule
    \end{tabular}
\quad 
    \begin{tabular}{c|c|cc}
        \multicolumn{4}{c}{$D=1$, quadratic} \\
         \midrule
         & iPCE & \multicolumn{2}{c}{niPCE} \\
         & $M$ & $M$ & $q$ \\
         \midrule
         \textbf{EE} & 10000 & 20000 & 50 \\
         \textbf{ETD-RDP} & 400 & 200 & 50\\
         \textbf{ETDRK4} & 200 & 100 & 50\\
         \textbf{ETDRK4} ref. & & 1000 & 200  \\
         \bottomrule
    \end{tabular}\\
    \begin{tabular}{c|c|cc}
        \multicolumn{4}{c}{$D=0$, cubic} \\
         \midrule
         & iPCE & \multicolumn{2}{c}{niPCE} \\
         & $M$ & $M$ & $q$ \\
         \midrule
         \textbf{EE} & 1000 & 500 & 50\\
         \textbf{ETD-RDP} & 200 & 200 & 50\\
         \textbf{ETDRK4} & 100 & 100 & 50\\
         \textbf{ETDRK4} ref. & & 1000 & 200  \\
         \bottomrule
    \end{tabular}
\quad
    \begin{tabular}{c|c|cc}
        \multicolumn{4}{c}{$D=1$, cubic} \\
         \midrule
         & iPCE & \multicolumn{2}{c}{niPCE} \\
         & $M$ & $M$ & $q$ \\
         \midrule
         \textbf{EE} & 20000 & 20000 & 50\\
         \textbf{ETD-RDP} & 400 & 200 & 50\\
         \textbf{ETDRK4} & 200 & 100 & 50 \\
         \textbf{ETDRK4} ref. & & 1000 & 200  \\
         \bottomrule
    \end{tabular}
    \caption{Numbers of step sizes $M$ and of samples (for niPCE) for each simulation shown in Figures \ref{linear_D=0} to \ref{cubic_D=1}. `ETDRK4 ref.' refers to the reference solution used for that simulation. For the linear equation with $D=0$, the exact solution is known.  }\label{parametertable}
\end{table}
\printbibliography
\end{document}